\documentclass[smallextended]{svjour3}
\smartqed  
\usepackage{graphicx}
\usepackage{amsfonts}
\usepackage{epsfig}
\usepackage{amssymb}
\usepackage{amsmath}
\newtheorem{stheorem}{Theorem}
\numberwithin{stheorem}{section}
\newtheorem{sremark}{Remark}
\numberwithin{sremark}{section}
\def\Frac#1#2{\frac{\displaystyle{#1}}{\displaystyle{#2}}}
\numberwithin{equation}{section}

\begin{document}

\title{Computation of asymptotic expansions of turning point
problems via Cauchy's integral formula: Bessel functions
\thanks{The authors acknowledge support from {\emph{Ministerio de Econom\'{\i}a y
Competitividad}}, project MTM2015-67142-P (MINECO/FEDER, UE). A.G. and J.S. acknowledge support 
from {\emph{Ministerio de Econom\'{\i}a y Competitividad}}, project MTM2012-34787.
A.G. acknowledges the Fulbright/MEC Program for support during her stay at SDSU.
J.S.  acknowledges the Salvador de Madariaga Program for support during his stay at SDSU.
}}

\titlerunning{Computation of asymptotic solutions of turning point
problems: Bessel functions}    

\author{T. M. Dunster         \and
        A. Gil \and
        J. Segura
}

\authorrunning{T.M. Dunster, A. Gil, J. Segura}

\institute{T. M. Dunster \at
              Department of Mathematics and Statistics\\
              San Diego State University. 5500 Campanile Drive San Diego, CA, USA.\\
              \email{mdunster@mail.sdsu.edu}   
           \and
           A. Gil \at
              Departamento de Matem\'atica Aplicada y CC. de la Computaci\'on.\\
ETSI Caminos. Universidad de Cantabria. 39005-Santander, Spain.
           \at
           \emph{and}
            \at
             Department of Mathematics and Statistics\\
              San Diego State University. 5500 Campanile Drive San Diego, CA, USA.\\
              \email{amparo.gil@unican.es} 
\and
           J. Segura \at
              Departamento de Matem\'aticas, Estad\'{\i}stica y Computaci\'on.\\
            Universidad de Cantabria. 39005-Santander, Spain.
           \at
           \emph{and}
            \at
             Department of Mathematics and Statistics\\
              San Diego State University. 5500 Campanile Drive San Diego, CA, USA.\\
              \email{segurajj@unican.es}          
}

\date{Received: date / Accepted: date}

\maketitle

\begin{abstract}

Linear second order differential equations having a large real parameter and
turning point in the complex plane are considered. Classical asymptotic
expansions for solutions involve the Airy function and its derivative, along
with two infinite series, the coefficients of which are usually difficult to
compute. By considering the series as asymptotic expansions for two
explicitly defined analytic functions, Cauchy's integral formula is employed
to compute the coefficient functions to high order of accuracy. The method
employs a certain exponential form of Liouville-Green expansions for
solutions of the differential equation, as well as for the Airy function. We
illustrate the use of the method with the high accuracy computation of
Airy-type expansions of Bessel functions of complex argument.

\keywords{Turning point problems \and Asymptotic expansions \and  Bessel functions \and Numerical computation}
 \subclass{MSC 34E05 \and 34E20 \and 33C10 \and 33F05}
\end{abstract}

\institute{T. M. Dunster \at
              Department of Mathematics and Statistics\\
              San Diego State University. 5500 Campanile Drive San Diego, CA, USA.\\
              \email{mdunster@mail.sdsu.edu}                      \and
           A. Gil \at
              Departamento de Matem\'atica Aplicada y CC. de la Computaci\'on.\\
ETSI Caminos. Universidad de Cantabria. 39005-Santander, Spain.
           \at
           \emph{and}
            \at
             Department of Mathematics and Statistics\\
              San Diego State University. 5500 Campanile Drive San Diego, CA, USA.\\
              \email{amparo.gil@unican.es} 
\and
           J. Segura \at
              Departamento de Matem\'aticas, Estad\'{\i}stica y Computaci\'on.\\
           ETSI Caminos. Universidad de Cantabria. 39005-Santander, Spain.
           \at
           \emph{and}
            \at
             Department of Mathematics and Statistics\\
              San Diego State University. 5500 Campanile Drive San Diego, CA, USA.\\
              \email{segurajj@unican.es}          
}

\section{Introduction}

In this paper we study linear second order differential equations having a
simple turning point. Specifically, we consider the differential equation%
\begin{equation}
d^{2}w/dz^{2}=\left\{ {u^{2}f(z)+g(z)}\right\} w,  \label{eq0}
\end{equation}
where $u$ is positive and large, $f(z)$ has a simple zero (turning point) at 
$z=z_{0}$ (say), and $f(z)$ and $g(z)$ are analytic in an unbounded domain
containing the turning point.

This is a classical problem, with applications to numerous special
functions. To obtain asymptotic solutions, the Liouville transformation 
\begin{equation}
\frac{2}{3}\zeta ^{3/2}=\pm \int_{z_{0}}^{z}f^{1/2}(t)dt,\quad W=\zeta
^{-1/4}f^{1/4}(z)w,  \label{eq2}
\end{equation}
is applied, where either sign in front of the integral can be chosen. As a
result we transform (\ref{eq0}) to the form 
\begin{equation}
d^{2}W/d\zeta ^{2}=\left\{ u^{2}\zeta +\psi (\zeta )\right\} W,  \label{eq3}
\end{equation}
where 
\begin{equation}
\psi (\zeta )=\frac{5}{16\zeta ^{2}}+\left\{ 4f(z){f}^{\prime \prime }(z)-5{f%
}^{\prime 2}(z)\right\} \frac{\zeta }{16f^{3}(z)}+\frac{\zeta g(z)}{f(z)}.
\label{eq4}
\end{equation}
The lower integration limit in (\ref{eq2}) ensures that the turning point 
$z=z_{0}$ of (\ref{eq0}) is mapped to the turning point ${\zeta }=0$ of 
(\ref{eq3}). Throughout this paper we shall assume that this turning point is
bounded away from any other turning points or singularities of (\ref{eq0}),
equivalently $\psi (\zeta )$ is analytic for $0\leq \left\vert \zeta
\right\vert <R$ for some positive $R$ which is independent of $u$.

When the turning point $z_0$ is real and $f(z)$ is real on a real interval
around $z_0$, the sign in (\ref{eq2}) is usually chosen in such a way that
the new variable $\zeta$ is real when $z$ is real and in a neighborhood of
the turning point.

\begin{figure}[tbp]
\begin{center}
\epsfxsize=6cm \epsfbox{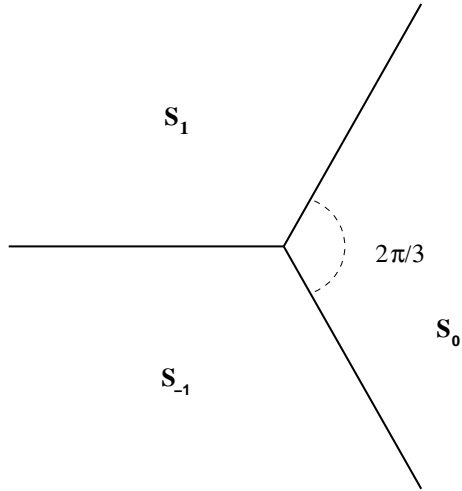}
\caption{The sectors ${\rm \bf S}_j$ in the complex plane}
\end{center}
\label{sectors}
\end{figure}

From \cite[Chap. 11, Theorem 9.1]{Olver:1997:ASF} we obtain solutions having
the following asymptotic expansions in terms of Airy functions 
\begin{equation}
\begin{array}{ll}
W_{2n+1,j}(u,\zeta )= & \mathrm{Ai}_{j}\left( u^{2/3}\zeta \right)
\displaystyle\sum_{s=0}^{n}\Frac{A_{s}(\zeta )}{u^{2s}} \\ 
& +\Frac{\mathrm{Ai}_{j}^{\prime }\left( u^{2/3}\zeta \right)}{u^{4/3}}
\displaystyle\sum\limits_{s=0}^{n-1}{\Frac{B_{s}(\zeta )}{u^{2s}}}+
\varepsilon_{2n+1,j}(u,\zeta ),
\end{array}
\label{eq5}
\end{equation}
for $j=0,\pm 1$. Here $\mathrm{Ai}_{j}(u^{2/3}\zeta )=\mathrm{Ai}%
(u^{2/3}\zeta e^{-2\pi ij/3})$, which are the Airy functions that are
recessive in the sectors $\mathrm{\mathbf{S}}_{j}:=\left\{ \zeta :|\mathrm{%
arg}(\zeta e^{-2\pi ij/3})|\leq \pi /3\right\} $ (Fig. \ref{sectors}); 
see \cite[\S 9.2(iii)]{Olver:2010:AAR}. Also, note that 
$\mathrm{Ai}_{j}^{\prime }(z)=d\mathrm{Ai}_{j}(z)/dz=e^{-2\pi ij/3}\mathrm{Ai}^{\prime }(ze^{-2\pi ij/3})
$.

From Olver's explicit error bounds we have 
$$
\varepsilon _{2n+1,j}(u,\zeta )=\mathrm{env}\left\{ \mathrm{Ai}_{j}\left(
u^{2/3}\zeta \right) \right\} \mathcal{O}\left( u^{-2n-1}\right) ,
$$
as $u\rightarrow \infty $, for $\zeta $ lying in certain domains described
in \cite[Chap. 11, \S 9]{Olver:1997:ASF}, and which we assume to be
unbounded. For a definition of the envelope function $\mathrm{env}$ for Airy
functions, see \cite[\S 2.8(iii)]{Olver:2010:AAR}.

In (\ref{eq5}) $A_{0}(\zeta)$ is an arbitrary non-zero constant (typically
taken to be 1), and for $s=0,1,2,\cdots $, the other coefficients satisfy
the recursion relations

\begin{equation}  \label{eq6a}
B_{s}(\zeta) =\frac{1}{2\zeta ^{1/2}}\int_{0}^{\zeta } \left\{\psi(t)
A_{s}(t) -A_{s}^{\prime \prime }(t) \right\} \frac{dt}{t^{1/2}},
\end{equation}
and 
\begin{equation}  \label{eq6b}
A_{s+1}(\zeta) =-\frac{1}{2}B_{s}^{\prime }(\zeta) +\frac{1}{2}\int
\psi(\zeta) B_{s}(\zeta) d\zeta .
\end{equation}

We remark that the lower integration limit in (\ref{eq6a}) must be $0$ in
order for each $B_{s}(\zeta) $ to be analytic at $\zeta =0$,
whereas in (\ref{eq6b}) there is no restriction in the choice of integration
constant. This will be of significance to us below.

In general, these coefficients are difficult to compute, primarily due to
the requirement of repeated integrations. They also show cancellations near
the turning point. For complex $\zeta $ close to $0$ one can compute these
coefficients in a numerically stable way by considering power series
expansions for the coefficients, as done in \cite{Amos:1986:A6A,Amos:1977:C6S}, 
where they are
expanded in powers of $\omega=\sqrt{1-z^2}$. For
other computational approaches to compute the coefficients, and in
particular for real values of $\zeta $, see \cite{Temme:1997:NAF}.

The purpose of this paper is to provide a more simple means of computing a
large number of these terms. We shall employ Cauchy's integral formula to do
so, and our results will be valid for real and complex $\zeta $ lying in a
bounded (but not necessarily small) domain containing the turning point 
$\zeta =0$. Our approach can potentially be extended to other situations,
including the cases of simple poles \cite[Chap. 12]{Olver:1997:ASF}, and
coalescing turning points.

We illustrate the use of the method with the high accuracy computation of
Airy-type expansions of Bessel functions of complex argument.

\section{General method}

We first present Liouville-Green expansions for solutions of (\ref{eq0}), a
certain form of which will be required for our method. These only involve
elementary (exponential) functions, but are not valid at the turning point.
The appropriate Liouville-Green transformation is given by \cite[Chap. 10, 
\S 2]{Olver:1997:ASF}, namely

\begin{equation}
\xi =\frac{2}{3}\zeta ^{3/2},\ V=f^{1/4}(z)w.  \label{eq7}
\end{equation}%
With these, equation (\ref{eq0}) is transformed to equation 
\begin{equation}
d^{2}V/d\xi ^{2}=\left\{ u^{2}+\phi (\xi )\right\} V,  
\end{equation}
where 
\begin{equation}
\phi (\xi )=\frac{4f(z)f^{\prime \prime }(z)-5f^{\prime 2}(z)}{16f^{3}(z)}+%
\frac{g(z)}{f(z)}.  \label{eq9}
\end{equation}

The branch for the first of (\ref{eq7}) will be dependent on the solutions
under consideration, as described below. Note, as $\zeta $ completes one
circuit about the turning point $\zeta =0$, the variable $\xi $
correspondingly crosses more than one Riemann sheet.

It turns out that solutions where asymptotic expansions appear inside
exponentials are more convenient for our purposes. Specifically, from \cite[%
Chap. 10, Ex. 2.1]{Olver:1997:ASF} we have solutions

\begin{equation}
V_{n}^{\pm }(u,\xi )=\exp \left\{ \pm u\xi +\displaystyle\sum_{s=1}^{n-1}
(\pm 1)^{s}\frac{E_{s}(\xi )}{u^{s}}\right\} +
\varepsilon_{n}^{\pm }(u,\xi ).
\label{eq10}
\end{equation}
In these, the coefficients are given by 
\begin{equation}
E_{s}(\xi )=\int F_{s}(\xi )d\xi \quad (s=1,2,3,\cdots ),  \label{eq12}
\end{equation}
where 
\begin{equation}
F_{1}(\xi )=\frac{1}{2}\phi (\xi ),\quad F_{2}(\xi )=-\frac{1}{4}\phi
^{\prime }(\xi ),  \label{eq13}
\end{equation}
and 
\begin{equation}
F_{s+1}(\xi )=-\frac{1}{2}F_{s}^{\prime }(\xi )-\frac{1}{2}
\displaystyle\sum_{j=1}^{s-1}F_{j}(\xi )F_{s-j}(\xi )\quad (s=2,3,\cdots ).  
\label{eq14}
\end{equation}
Primes are derivatives with respect to $\xi $. The integration constants in 
(\ref{eq12}) will be discussed below, and we find that for the odd
coefficients $E_{2j+1}(\xi )$ ($j=0,1,2,\cdots $) they must be suitably
chosen.

Explicit error bounds for $\varepsilon_{n}^{\pm }(u,\xi )$, which verify
the asymptotic validity of the expansions (\ref{eq10}), are given in 
\cite{Dunster:1998:AOT}. In particular, for arbitrary $\delta >0$, under certain
conditions on $\psi (\xi )$, we have that $\varepsilon_{n}^{\pm }(u,\xi)
=e^{\pm u\xi }\mathcal{O}\left( u^{-n}\right) $ as $\xi \rightarrow \infty$
in certain domains $\Xi ^{\pm }$ as described in 
\cite[Chap. 10, \S 3]{Olver:1997:ASF}. These are the same as those for the corresponding asymptotic solutions of the more common form 
\begin{equation}
V\sim \exp \{\pm u\xi \}\sum_{s=0}^{\infty }(\pm 1)^{s}A_{s}(\xi
)u^{-s}.  \label{eq11a}
\end{equation}

It is the relation (\ref{eq14}) that is the reason why the expansions 
(\ref{eq10}) are numerically advantageous: the coefficients $F_{s}(\xi )$ 
can all be determined explicitly without resorting to
integration. Furthermore, from (\ref{eq12}) we observe that only one
integration is required (numerical or explicit) to evaluate each 
$E_{s}(\xi) $, as opposed to repeated integrals for computing the coefficients 
$A_{s}(\xi )$ in (\ref{eq11a}).

Remarkably, it turns out that integration is not required to evaluate the
even terms $E_{2j}(\xi )$ ($j=1,2,3,\cdots $). To see this, consider the
Wronskian of the solutions $V_{n}^{\pm }(u,\xi )$ given by (\ref{eq10}).
Since this is a constant (by Abel's theorem) we infer that
$$
\left\{ u+\sum_{j=0}^{\infty }\frac{F_{2j+1}(\xi )}{u^{2j+1}}\right\} \exp
\left\{ 2\sum_{j=1}^{\infty }\frac{E_{2j}(\xi )}{u^{2j}}\right\} \sim
\mbox{constant,}
$$
which, on taking logarithms, yields 
\begin{equation}
\sum_{j=1}^{\infty }\frac{E_{2j}(\xi )}{u^{2j}}\sim 
-\frac{1}{2}\ln \left\{1+\sum_{j=0}^{\infty }
\frac{F_{2j+1}(\xi )}{u^{2j+2}}\right\} +
\mbox{constant.}  
\label{wronsk}
\end{equation}
We then asymptotically expand the RHS of this relation in inverse powers of 
$u^{2}$, and equate the coefficients of both sides. As a result, we find that
(to within an arbitrary additive constant in each instance) 
\begin{equation}
E_{2}(\xi )=-\frac{1}{2}F_{1}(\xi ),\,E_{4}(\xi )=\frac{1}{4}F_{1}^{2}(\xi )
-\frac{1}{2}F_{3}(\xi ),  
\end{equation}
and so on. In particular, the even coefficients $E_{2j}(\xi )$ 
($j=1,2,3,\cdots $) are explicitly given in terms of $F_{2k+1}(\xi )$ 
($k=0,1,2,\cdots ,j-1$) (which in turn are given by (\ref{eq13}) and 
(\ref{eq14})).

At this stage we consider Liouville-Green solutions of (\ref{eq3}).
Comparing (\ref{eq2}), (\ref{eq7}) and (\ref{eq10}) we obtain three
asymptotic solutions $W_{j}(u,\zeta )$ ($j=0,\pm 1$) of (\ref{eq3}) which
are recessive (respectively) for $\zeta \in \mathrm{\mathbf{S}}_{j}$,
possessing the asymptotic expansions 
\begin{equation}
W_{0}(u,\zeta )\sim {\zeta ^{-1/4}}\exp \left\{ -\frac{2}{3}u\zeta ^{3/2}+
\displaystyle\sum_{s=1}^{\infty }(-1)^{s}\frac{E_{s}(\xi )}{u^{s}}\right\} ,
\label{eq16}
\end{equation}
and
\begin{equation}
W_{j}(u,\zeta )\sim {\zeta ^{-1/4}}\exp \left\{ \frac{2}{3}u\zeta ^{3/2}+%
\displaystyle\sum_{s=1}^{\infty }\frac{E_{s}(\xi )}{u^{s}}\right\} \ (j=\pm
1),  \label{eq16a}
\end{equation}
as $u\rightarrow \infty ,$ uniformly for $\zeta \in \Omega _{j}$\ (say). For 
$j=0,\pm 1$ the branches are such that ${\rm Re}\left( e^{-ij\pi }\zeta
^{3/2}\right) \geq 0$ for $\zeta \in \mathrm{\mathbf{S}}_{j}$ , and 
${\rm Re}\left( e^{-ij\pi }\zeta ^{3/2}\right) \leq 0$ for $\zeta \notin 
\mathrm{\mathbf{S}}_{j}$. As is typically the case in practice, for each $j$ we
assume that $\Omega _{j}\cap \mathrm{\mathbf{S}}_{j}$ is unbounded.

We remark that, on account of the analyticity of $\psi (\zeta )$ in the disk 
$0\leq \left\vert \zeta \right\vert <R$, the expansions (\ref{eq16}) and 
(\ref{eq16a}) certainly hold in the bounded sector $\delta \leq \left\vert \zeta
\right\vert <R$, $|\mathrm{arg}(\zeta e^{-2\pi ij/3})|\leq \pi -\delta $,
where here and elsewhere $\delta $ denotes an arbitrary small positive
constant. We also note that the recessive property at $\zeta
=\infty $ in $\Omega _{j}\cap \mathrm{\mathbf{S}}_{j}$ uniquely defines 
$W_{j}(u,\zeta )$ up to a multiplicative constant. Indeed, we have that 
$W_{j}(u,\zeta )=c_{2n+1,j}W_{2n+1,j}(u,\zeta )$ for some constants 
$c_{2n+1,j}$, although we shall not use these relations.

Now, since no two from these three solutions are linearly dependent, we can
assume they satisfy a connection formula of the form
\begin{equation}
\lambda _{-1}W_{-1}(u,\zeta )=iW_{0}(u,\zeta )+\lambda _{1}W_{1}(u,\zeta ),
\label{eq21}
\end{equation}
for certain constants $\lambda _{-1}$ and $\lambda _{1}$ (which may of
course depend on $u$). The factor $i$ is for convenience.

We note that each $W_{j}(u,\zeta )$ ($j=0,\pm 1$), being a solution of 
(\ref{eq0}), is analytic in a neighborhood of the turning point. Based on 
(\ref{eq5}), and following \cite{Boyd:1987:AEF,Temme:1997:NAF}, we thus can 
\textit{define} functions $A(u,z)$ and $B(u,z)$, analytic at $z=z_{0}$ 
($\zeta =0$), implicitly by the pair of equations 
\begin{equation}
\Frac{1}{2\pi ^{1/2}u^{1/6}}W_{0}(u,\zeta )=
\mathrm{Ai}_{0}(u^{2/3}\zeta )A(u,z)+
\mathrm{Ai}_{0}^{\prime }(u^{2/3}\zeta)B(u,z),  \label{eq17}
\end{equation}
and
\begin{equation}
\frac{e^{\pi i/6}\lambda_1}{2\pi^{1/2}u^{1/6}}
W_{1}(u,\zeta )=\mathrm{Ai}_{1}(u^{2/3}\zeta )A(u,z)+
\mathrm{Ai}_{1}^{\prime }(u^{2/3}\zeta )B(u,z),  \label{eq17a}
\end{equation}
where (as shown below) the multiplicative constants on the LHS of both
equations have been chosen to yield the appropriate behavior of $A(u,z)$ and 
$B(u,z)$ as $u\rightarrow \infty $. We remark that for computational
purposes it is more convenient to consider $A(u,z)$ and $B(u,z)$ as
functions of $z$, although in deriving their asymptotic expansions we shall
regard them as functions of $\zeta $ as necessary.

Next, from the connection formula (\ref{eq21}), and the corresponding
well-known connection formula for the Airy functions
$$
\mathrm{Ai}_{-1}\left(u^{2/3}\zeta \right) =e^{-\pi i/3}
\mathrm{Ai}_{1}\left( u^{2/3}\zeta \right) 
+e^{\pi i/3}\mathrm{Ai}_{0}\left(u^{2/3}\zeta \right) ,  
$$
we derive from (\ref{eq17}) and (\ref{eq17a}) the following Airy function
representation for the solution which is recessive in 
$\mathrm{\mathbf{S}}_{-1}$ 
\begin{equation}
\Frac{e^{-\pi i/6}\lambda_{-1}}{2\pi^{1/2}u^{1/6}}
W_{-1}(u,\zeta )=\mathrm{Ai}_{-1}(u^{2/3}\zeta )A(u,z)+
\mathrm{Ai}_{-1}^{\prime }(u^{2/3}\zeta )B(u,z).  \label{eq17b}
\end{equation}

We shall show, using (\ref{eq17}), (\ref{eq17a}) and (\ref{eq17b}), that 
$A(u,z)$ and $B(u,z)$ are slowly varying relative to the fast variation of
the Airy functions in a \textit{full} neighborhood of turning point.
Specifically, on referring to (\ref{eq5}), $A(u,z)$ and $B(u,z)$ will admit
the following asymptotic expansions as $u\rightarrow \infty $ 
\begin{equation}
A(u,z)\sim \displaystyle\sum_{s=0}^{\infty }\frac{A_{s}(\zeta) 
}{u^{2s}},\,B(u,z)\sim \frac{1}{u^{4/3}}\displaystyle\sum_{s=0}^{\infty }%
\frac{B_{s}(\zeta) }{u^{2s}},  \label{expannu}
\end{equation}
uniformly with respect to $\zeta $ lying in a certain unbounded domain which
contains the disk $0\leq \left\vert \zeta \right\vert <R$.

To show the slowly-varying nature of $A(u,z)$ and $B(u,z)$ (for example)
when $\zeta \in \left\{ \mathrm{\mathbf{S}}_{0}\cap \Omega _{0}\right\} \cup
\left\{ \mathrm{\mathbf{S}}_{-1}\cap \Omega _{-1}\right\} $, we solve (\ref%
{eq17}) and (\ref{eq17b}) for the coefficients, and using the Wronskian for
Airy functions \cite[\S 9.2(iv)]{Olver:2010:AAR}, we arrive at the explicit
representations 
\begin{equation}
\begin{array}{ll}
A(u,z)=\Frac{\pi ^{1/2}}{u^{1/6}} & \left\{ e^{\pi i/6}W_{0}(u,\zeta )
\mathrm{Ai}_{-1}^{\prime }(u^{2/3}\zeta )\right. \\ 
& \left. -\lambda _{-1}W_{-1}(u,\zeta )
\mathrm{Ai}_{0}^{\prime}(u^{2/3}\zeta )\right\} ,
\end{array}
\label{eq18}
\end{equation}
and
\begin{equation}
\begin{array}{ll}
B(u,z)=\Frac{\pi ^{1/2}}{u^{1/6}} & \left\{ \lambda _{-1}W_{-1}(u,\zeta )
\mathrm{Ai}_{0}(u^{2/3}\zeta )\right. \\ 
& \left. -e^{\pi i/6}W_{0}(u,\zeta )\mathrm{Ai}_{-1}(u^{2/3}\zeta )\right\} .
\end{array}
\label{eq19a}
\end{equation}
Then for $\zeta \in \left\{ \mathrm{\mathbf{S}}_{0}\cap \Omega _{0}\right\}
\cup \left\{ \mathrm{\mathbf{S}}_{-1}\cap \Omega _{-1}\right\} $ each
product pair of functions on the RHS of (\ref{eq18}) and (\ref{eq19a})
consists of one exponentially small function times an exponentially large
one. Consequently these forms are numerically satisfactory in that domain.

At this stage the integration constants in (\ref{eq12}) are arbitrary, and
indeed do not have to be independent of $u$. Therefore, for notational
convenience, in (\ref{eq18}) and (\ref{eq19a}) we shall absorb the
coefficient $\lambda _{-1}$ into the expansion\ (\ref{eq16a}) for 
$W_{-1}(u,\zeta )$, by redefining the Liouville-Green coefficients if
necessary. For example,we can redefine $E_{1}(\xi )$ to be 
$E_{1}(\xi )-\frac{1}{2}u\ln \left( \lambda _{-1}\right) $ and 
$E_{2}(\xi )$ to be $E_{2}(\xi )-\frac{1}{2}u^{2}\ln \left( \lambda _{-1}\right)$: 
as a result $\lambda _{-1}W_{-1}(u,\zeta )$ is scaled to become $W_{-1}(u,\zeta )$,
whereas from (\ref{eq16}) we see that $W_{0}(u,\zeta )$ is unchanged.

We now use (\ref{A9}) and (\ref{A12}), along with the corresponding
expansions (which are valid for $\left\vert {\arg \left( {\zeta e^{2\pi i/3}}%
\right) }\right\vert \leq \pi -\delta $) 
$$
\mathrm{Ai}_{-1}\left(u^{2/3}\zeta\right) \sim \frac{e^{-\pi i/6}}{2\pi
^{1/2}u^{1/6}\zeta ^{1/4}}\exp \left\{ u\xi +\sum\limits_{s=1}^{\infty }
\frac{a_{s}}{su^{s}\xi ^{s}}\right\} ,  
$$
and 
$$
\mathrm{Ai}_{-1}^{\prime }\left(u^{2/3}\zeta\right) 
\sim \frac{e^{-\pi i/6}u^{1/6}\zeta ^{1/4}}{2\pi ^{1/2}}\exp 
\left\{u\xi +\sum\limits_{s=1}^{\infty }\frac{\tilde{{a}}_{s}}{su^{s}\xi ^{s}}
\right\} .  
$$
Then, from (\ref{eq18}) and (\ref{eq19a}), along with the defining
expansions (\ref{eq10}), we obtain our main result, the slowly-varying expansions

\begin{equation}
\begin{array}{ll}
A(u,z)\sim & \exp \left\{ \displaystyle\sum_{j=1}^{\infty }
\Frac{E_{2j}(\xi )+\tilde{a}_{2j}\xi ^{-2j}/(2j)}{u^{2j}}\right\} \\ 
& \times \mathrm{cosh}\left\{ \displaystyle\sum_{j=0}^{\infty }
\Frac{E_{2j+1}(\xi )-\tilde{a}_{2j+1}\xi ^{-2j-1}/(2j+1)}
{u^{2j+1}}\right\} ,
\end{array}
\label{eq29}
\end{equation}
and

\begin{equation}
\begin{array}{ll}
B(u,z)\sim & \Frac{1}{u^{1/3}\zeta ^{1/2}}
\exp\left\{ \displaystyle\sum_{j=1}^{\infty }
\Frac{E_{2j}(\xi)+a_{2j}\xi ^{-2j}/(2j)}{u^{2j}}\right\} \\ 
& \times \mathrm{sinh}\left\{ \displaystyle\sum_{j=0}^{\infty }
\Frac{E_{2j+1}(\xi )-a_{2j+1}\xi ^{-2j-1}/(2j+1)}{u^{2j+1}}\right\} .
\end{array}
\label{eq30}
\end{equation}

These asymptotic expansions certainly hold for (at least) $-\pi +\delta \leq 
\mathrm{arg}(\zeta )\leq \frac{1}{3}\pi -\delta ,$ $0<\left\vert \zeta
\right\vert <R$. Similar slowing-varying expansions can be obtained for all
other values of $\mathrm{arg}(\zeta )$ near the turning point, by solving
for $A(u,z)$ and $B(u,z)$ from a suitably chosen pair from the three
equations (\ref{eq17}), (\ref{eq17a}) and (\ref{eq17b}). Each of these
expansions can be expressed in the same form as above, except that the
coefficients $E_{s}(\xi )$ may differ by the integration constants in 
(\ref{eq12}).

In fact, if these integration constants are arbitrarily chosen, we find that
the RHS of (\ref{eq29}) and (\ref{eq30}) generally each has a branch point
at $\zeta =0$ ($z=z_{0}$). Hence, being multi-valued (in fact unbounded) for
small $\zeta $, the expansions are generally only valid for restricted 
$\arg(\zeta)$ near the turning point. It is only for specific
integration constants in (\ref{eq12}), at least for the odd coefficients 
$E_{2j+1}(\xi )$ ($j=0,1,2\cdots $), that the expansions are single-valued,
and hence by a continuity argument (\ref{eq29}) and (\ref{eq30}) hold for 
$0<\left\vert \zeta \right\vert <R$, with $\mathrm{arg}(\zeta )$
unrestricted. This essentially corresponds to the choice of the lower
integration constant in (\ref{eq6a}), which (as remarked above) ensures that
each coefficient $B_{s}(\zeta )$ is analytic at $\zeta =0$.

Recalling that $\xi =\frac{2}{3}\zeta ^{3/2}$, it is straightforward to show that
the expansions (\ref{eq29}) and (\ref{eq30}) are single-valued at $\zeta =0$
if $\zeta ^{1/2}E_{2j+1}(\xi )$ and $E_{2j}(\xi )$ ($j=0,1,2,\cdots $) are
meromorphic: here and throughout meromorphic means with respect to $\zeta $,
and at $\zeta =0$. If we assume for the moment that this is true for 
$\zeta^{1/2}E_{2j+1}(\xi )$ and $E_{2j}(\xi )$, then on re-expanding (\ref{eq29})
and (\ref{eq30}) into the forms (\ref{expannu}), we deduce that the
coefficients $A_{s}(\zeta) $ and $B_{s}(\zeta)$
in the latter expansions must be single-valued. By a uniqueness argument
these coefficients must satisfy (\ref{eq6a}) and  (\ref{eq6b}), and moreover
single-valuedness means the lower integration limit in the former must
indeed be 0. But in this case we know that  $A_{s}(\zeta)$ and 
$B_{s}(\zeta)$ are actually analytic at $\zeta =0$. We deduce
that if $\zeta ^{1/2}E_{2j+1}(\xi )$ and $E_{2j}(\xi )$ are meromorphic then
(\ref{eq29}) and (\ref{eq30}) have a removable singularity when they are
re-expanded into the forms (\ref{expannu}). With an appropriate choice of
integration constant for each $E_{2j+1}(\xi )$ ($j=0,1,2,\cdots $) we now
show that this is indeed so.

To this end, we first note from (\ref{eq4}) and (\ref{eq9}) that 
\begin{equation}
\phi (\xi )=\frac{\psi (\zeta )}{\zeta }-\frac{5}{16\zeta ^{3}},
\end{equation}
and hence this function is meromorphic. Therefore, on noting that 
$d\xi/d\zeta =\zeta ^{1/2}$, we find by induction from 
(\ref{eq12}) - (\ref{eq14}) that $\zeta ^{1/2}F_{2j}(\xi )$ and 
$F_{2j+1}(\xi )$ are also meromorphic.

To establish the desired meromorphicity of $\zeta ^{1/2}E_{2j+1}(\xi )$ and 
$E_{2j}(\xi )$, let us consider these separately. Firstly, from (\ref{eq12})
we see that $\zeta ^{1/2}{E_{2j+1}(\xi )}$ is meromorphic if and only if the
Laurent-type expansion 
\begin{equation}
\zeta ^{1/2}E_{2j+1}(\xi )=\zeta ^{1/2}\int \zeta ^{1/2}F_{2j+1}(\xi
(\zeta )) d\zeta =\alpha _{2j+1}\zeta ^{1/2}+\sum_{k=-2j-1}^{\infty
}\alpha _{2j+1,k}\zeta ^{k},  \label{eq30a}
\end{equation}
satisfies $\alpha _{2j+1}=0$. Clearly the integration constants in 
(\ref{eq12}) can be selected in order for this to be so,\textit{ and we assume
this from now on}.

Next consider the even terms. Again from (\ref{eq12}), we observe that 
\begin{equation}
E_{2j}(\xi )=\int \zeta ^{1/2}F_{2j}(\xi (\zeta )) d\zeta ,
\label{Ezeta}
\end{equation}
must either be meromorphic, or have a logarithmic singularity at $\zeta =0$.
However, the latter possibility can immediately be discarded by referring to
(\ref{wronsk}) along with the meromorphicity of $F_{2j+1}(\xi )$. We note
that the integration constants in (\ref{Ezeta}) do not affect the
meromorphicity of $E_{2j}(\xi )$, and hence they can be arbitrarily chosen.

In summary, we have established that $\alpha _{2j+1}=0$ in (\ref{eq30a}) is
a necessary and sufficient condition for the expansions (\ref{eq29}) and 
(\ref{eq30}) to be valid in a neighborhood of the turning point with 
$\mathrm{arg}(\zeta )$ unrestricted.

In practice there are several ways of ensuring the correct choice of $%
E_{2j+1}(\xi )$. If each of these functions can be explicitly determined
from (\ref{eq12}), then a symbolic algebra system could be used to determine
the expansion (\ref{eq30a}), and if $\alpha _{2j+1}\neq 0$ one can simply
subtract this constant from the function given by the anti-derivative
derived from (\ref{eq12}).

If explicit integration of (\ref{eq12}) is not possible, and quadrature is
required, then from (\ref{eq30a}) we observe that 
\begin{equation}
\alpha _{2j+1}=\frac{1}{2}\left\{ E_{2j+1}(\xi ^{\ast })+E_{2j+1}(\xi
)\right\} ,  \label{lambda}
\end{equation}%
where $\xi^\ast=\xi(\zeta e^{2\pi i})$. Since the
coefficients will be computed numerically on a loop surrounding $\zeta =0$,
the values on the RHS of (\ref{lambda}) can also be computed, provided $\xi $
and $\xi ^{\ast }$ correspond to the initial and terminal points of the
loop. Hence, as in the case of the explicitly known coefficients, these
numerically computed values can be subtracted from the initial values of 
$E_{2j+1}(\xi )$ evaluated from (\ref{eq12}).

In the common situation where $\psi (\zeta )$ is real in some real interval 
$(-a,a)$ where $a>0$, this numerical method can be simplified somewhat.
Specifically, if we choose the lower limit of integration in (\ref{eq12}) to
be at a real point on our path of integration, $\xi =\xi (\zeta _{0})$ say,
where $0<\zeta _{0}<a$, then we similarly find from (\ref{eq30a}) that
\begin{equation}
\alpha _{2j+1}=\mathrm{Re}\left\{ E_{2j+1}\left( \xi (\zeta _{0}e^{\pi
i})\right) \right\} ,  
\end{equation}
and we can then proceed as above.

We summarize the principal result as follows.

\begin{stheorem}
For the differential equation
\begin{equation}
d^{2}w/dz^{2}=\left\{ {u^{2}f(z)+g(z)}\right\} w,  
\label{thmode}
\end{equation}
assume $u$ is positive and large, $f(z)$ has a
simple zero at  $z=z_{0}$, and $f(z)$ and $g(z)$ 
are analytic in a domain $D$ containing $z_{0}$. 
Further assume that $f(z)$ does not vanish in the disk 
$D\left(z_{0},\rho \right) :=\left\{ z:0<\left\vert z-z_{0}\right\vert <\rho
\right\} \subset D$.
Define variables $\xi$ and $\zeta$ by
\begin{equation}
\xi =\frac{2}{3}\zeta ^{3/2}=\int_{z_{0}}^{z}f^{1/2}(t)dt,  
\label{xizeta}
\end{equation}
and let $\mathrm{Ai}_{j}(u^{2/3}\zeta )$ ($j=0,\pm 1$) 
denote the Airy functions $\mathrm{Ai}(u^{2/3}\zeta e^{-2\pi ij/3})$. 
Then there exist three numerically satisfactory solutions
of (\ref{thmode}) given by
\begin{equation}
w_{j}\left( u,z\right) =\zeta ^{1/4}f^{-1/4}(z)\left\{ \mathrm{Ai}%
_{j}(u^{2/3}\zeta )A(u,z)+\mathrm{Ai}_{j}^{\prime }(u^{2/3}\zeta
)B(u,z)\right\} .  
\end{equation}
In these, the coefficient functions $A(u,z)$ and $B(u,z)$
are analytic at $z=z_{0}$, and possess the asymptotic
expansions (\ref{eq29}) and (\ref{eq30}) in a domain that includes 
$D\left(z_{0},\rho \right)$. Here $a_{1}=a_{2}=\frac{5}{72}$,
$\tilde{a}_{1}=\tilde{a}_{2}=-\frac{7}{72}$, 
and in both cases subsequent terms are given by (\ref{A7}). 
The coefficients $E_{s}(\xi )$ ($s=1,2,3,\cdots$) are given by 
(\ref{eq12}) - (\ref{eq14}), where the integration constants for the odd 
coefficients in (\ref{eq12}) must be selected so that 
$\zeta ^{1/2}E_{2j+1}(\xi )$ ($j=0,1,2,\cdots $) 
is meromorphic as a function of $\zeta$
at $\zeta =0$; i.e. $\zeta ^{1/2}E_{2j+1}(\xi )$
possesses the Laurent expansion (\ref{eq30a}) with $\alpha _{2j+1}=0$.
\end{stheorem}

\begin{sremark}
The even terms $E_{2j}(\xi )$ ($j=1,2,3,\cdots $) can alternatively be
evaluated by expanding the RHS of (\ref{wronsk}) in inverse powers of $u^{2}$, 
and equating the coefficients of both sides.
This avoids integration for evaluating these terms. We also note that the
expansions (\ref{eq29}) and (\ref{eq30}) are valid in the same domains as
those for the corresponding asymptotic solutions of \cite[Chap. 11, Theorem
9.1]{Olver:1997:ASF}, and can be unbounded provided $f(z)$ and $g(z)$ have
the appropriate behavior at $\infty $ (see \cite[Chap. 11, Sect. 9.3]%
{Olver:1997:ASF}).
\end{sremark}

Our focus is to utilize the expansions (\ref{eq29}) and (\ref{eq30}) to
efficiently compute $A(u,z)$ and $u^{4/3}B(u,z)$ to 
$\mathcal{O}\left(u^{-2m}\right) $ for some prescribed $m$. As we shall 
see in the application
to Bessel functions below, in general it may be more convenient to consider
scaled functions 
\begin{equation}
\mathcal{A}(u,z)=\phi (u,z)A(u,z),\,\mathcal{B}(u,z)=\phi (u,z)B(u,z),
\label{scaled}
\end{equation}
where $\phi (u,z)$ is some suitably chosen function which is
analytic in a domain $\Omega $ (say) in which the expansion (\ref{eq29}) is
valid.

Our aim is to employ the Cauchy integrals 
\begin{equation}
\mathcal{A}(u,z)=\frac{1}{2\pi i}\oint\limits_{\mathcal{L}}\frac{\phi
(u,t)A(u,t)}{t-z}dt,\,\mathcal{B}(u,z)=\frac{1}{2\pi i}
\oint\limits_{\mathcal{L}}\frac{\phi (u,t)B(u,t)}{t-z}dt,  \label{eq31}
\end{equation}
where $\mathcal{L}$ is a positively orientated closed loop lying in $\Omega $
and surrounding $t=z$ and $t=z_{0}$, and therefore (\ref{eq29}) and 
(\ref{eq30}) can be inserted in the respective integrands. As we showed above,
the asymptotic expansions (\ref{eq29}) and (\ref{eq30}) are in theory valid
close to the turning point $\xi =\zeta =0$ ($z=z_{0}$) as $u\rightarrow
\infty $. However, for fixed $u$ the terms in the series become unbounded as 
$\xi \rightarrow 0$, and therefore these series are numerically
unsatisfactory near the turning point. It is for this reason that we use 
the Cauchy integral formulas (\ref{eq31}) for their numerical approximation.

At this point, it is important to realize that we will needed to compute the
coefficients $E_{s}(\xi )$ in (\ref{eq29}) and (\ref{eq30}) at specific
fixed points on $\mathcal{L}$. Since this computation is done once and for
all the computational efficiency in the computation of these coefficients is
not so important. And, as noted earlier, if the integrals in (\ref{eq12})
cannot be explicitly evaluated, only one numerical evaluation of an integral
is required for each odd coefficient.

For our purposes it will be sufficient to consider a circular path of
integration with center at $z_c$, not necessarily with $z_c=z_0$, but which in any
case must contain the turning point $z_0$. For computing the integral, we therefore use
the parametrization $t(\theta )=z_{c}+Re^{i\theta }$, $\theta \in \lbrack
0,2\pi ]$ and we have 
\begin{equation}
\mathcal{A}(u,z)=\frac{1}{2\pi }\int_{0}^{2\pi }F(\theta )d\theta
,\,F(\theta )=\mathcal{A}(u,t(\theta ))\frac{t(\theta )-z_{c}}{t(\theta )-z},
\label{31a}
\end{equation}%
and likewise for $\mathcal{B}(u,z)$. The function $F(\theta )$ is periodic
and we are integrating over one period. 
Since the function is infinitely differentiable as a function of $\theta $ we can expect that
the trapezoidal rule will give good convergence \cite[Thm. 5.6]{Gil:2007:NSF}. 
We write $\theta _{j}=2\pi j/N$ ( $j=0,1,\dots N$) and we have 
$$
\int_{0}^{2\pi }F(\theta )d\theta \approx \frac{2\pi }{N}\left\{ \frac{1}{2}%
(F(\theta _{0})+F(\theta _{N}))+\sum_{j=1}^{N-1}F(\theta _{j})\right\} ,
$$
and therefore 
\begin{equation}
\mathcal{A}(u,z)\approx \frac{1}{N}\sum_{j=1}^{N}F(2\pi j/N).  \label{31c}
\end{equation}

Numerical experiments show that $N$ does not need to be large and that, for
instance in our application to Bessel equation (see \S \ref{ATE}), $N=500$ is enough
for more than $15$-digits accuracy in a wide region around the turning point 
(namely for $z_c=2$, $R=1.8$); the number can be reduced for smaller regions, 
and for $z_c=1$, $R=0.5$, $15$ digits accuracy is reached with $N=150$. This is consistent with the expected good performance of the trapezoidal rule.

For the derivatives of solutions $w_{j}(u,z) $ (say) of 
(\ref{eq0}), suppose scaled coefficient functions $\mathcal{A}(u,z)$ 
and $\mathcal{B}(u,z) $ are defined by 
\begin{equation}
\lambda _{j}w_{j}(u,z)=\mathrm{Ai}_{j}\left( u^{2/3}\zeta \right) 
\mathcal{A}(u,z)+\mathrm{Ai}_{j}^{\prime }\left( u^{2/3}\zeta \right) 
\mathcal{B}(u,z)\ (j=\pm 1) ,  \label{31d}
\end{equation}%
for some appropriate connection coefficients $\lambda _{\pm 1}$. We can then
proceed similarly as before, by defining corresponding coefficients 
$\mathcal{C}(u,z) $ and $\mathcal{D}(u,z)$ by
\begin{equation}
\lambda _{j}w_{j}^{\prime }(u,z)=\mathrm{Ai}_{j}\left( u^{2/3}\zeta \right) 
\mathcal{C}(u,z)+\mathrm{Ai}_{j}^{\prime }\left( u^{2/3}\zeta \right) 
\mathcal{D}(u,z)\ (j=\pm 1) .  \label{31e}
\end{equation}%
Then, solving for these coefficients, we obtain the same expressions as for
the previous coefficients, but with $w_{j}\left( {u,z}\right) $ substituted
by their derivatives. Then, to obtain similar expansions to (\ref{eq29}) and
(\ref{eq30}), we would also need the Liouville-Green expansions for the
derivatives of these functions, which can be obtained by differentiation of
the corresponding expansions for the functions themselves.

A preferable way to compute the new coefficients is to differentiate 
(\ref{31d}). After differentiating with respect to $z$, using the Airy
differential equation to eliminate the second derivative of the Airy
functions, then solving and comparing to (\ref{31e}), we obtain the relations

\begin{equation}
\mathcal{C}(u,z)={\mathcal{A}}^{\prime }(u,z)+u^{4/3}\zeta (z)\zeta ^{\prime
}(z)\mathcal{B}(u,z),  \label{coeffc}
\end{equation}
and 
\begin{equation}
\mathcal{D}(u,z)=u^{2/3}\zeta ^{\prime }(z)\mathcal{A}(u,z)+ {\mathcal{B}}%
^{\prime }(u,z).  \label{coeffd}
\end{equation}

To compute the derivatives of the coefficients we can also use Cauchy's
integral formula to write 
\begin{equation}
\mathcal{A}^{\prime }(u,z)=\frac{1}{2\pi i}\oint\limits_{\mathcal{L}} 
\frac{\mathcal{A}(u,t)}{(t-z)^{2}}dt,\ \mathcal{B}^{\prime }(u,z)=
\frac{1}{2\pi i}\oint\limits_{\mathcal{L}} 
\frac{\mathcal{B}(u,t)}{(t-z)^{2}}dt .
\label{derica}
\end{equation}
Then, from (\ref{scaled}), (\ref{eq31}), and (\ref{coeffc}) - (\ref{derica})
we obtain 
\begin{equation}
\mathcal{C}(u,z)=\frac{1}{2\pi i}\oint\limits_{\mathcal{L}}
\frac{\mathcal{A}(u,t)+u^{4/3}{(t-z)}\zeta (t)\zeta ^{\prime }(t)
\mathcal{B}(u,t)}{(t-z)^{2}}dt,  \label{eq57}
\end{equation}%
and%
\begin{equation}
\mathcal{D}(u,z)=\frac{{1}}{{2\pi i}}\oint\limits_{\mathcal{L}}\frac{%
\mathcal{B}(u,t)+u^{2/3}{(t-z)}\zeta ^{\prime }(t)\mathcal{A}(u,t)}{{%
(t-z)^{2}}}dt.  \label{eq58}
\end{equation}%
We then use the computed values of $\mathcal{A}(u,t)$ and $\mathcal{B}(u,t)$
on $\mathcal{L}$, and proceed as above. Again, the trapezoidal rule is a
good choice; in fact, it is in some sense optimal \cite{Bornemann:2011:AAS}.

We remark that it is desirable that (\ref{coeffc}) and (\ref{coeffd}) be
numerically stable, in the sense that, for large $u$, there is no
cancellation in leading order terms in the two terms in either
representation. The choice of scaling function $\phi (u,z) $ in 
(\ref{scaled}) should be such that this is indeed the case.

\section{Bessel's equation: preliminary transformations}

We illustrate the new technique using the Airy function asymptotic
expansions for Bessel functions. The first step in doing so, is to apply the
Liouville transformations described in \S \S 1 and 2 to Bessel's equation.
To this end, we first note that functions $w=z^{1/2}J_{\nu }(\nu z)$, 
$w=z^{1/2}H_{\nu }^{(1)}(\nu z)$ and $w=z^{1/2}H_{\nu }^{(2)}(\nu z)$ satisfy 
$$
\frac{d^{2}w}{dz^{2}}=\left\{\nu ^{2}\frac{1-z^{2}}{z^{2}}-\frac{1}{4z^{2}}
\right\} w.  
$$
Here $\nu $ plays the role of our large positive parameter $u$, and $z$ is
complex.

\begin{figure}[tb]
\vspace*{0.8cm}
\par
\begin{center}
\begin{minipage}{3cm}
\centerline{\protect\hbox{\psfig{file=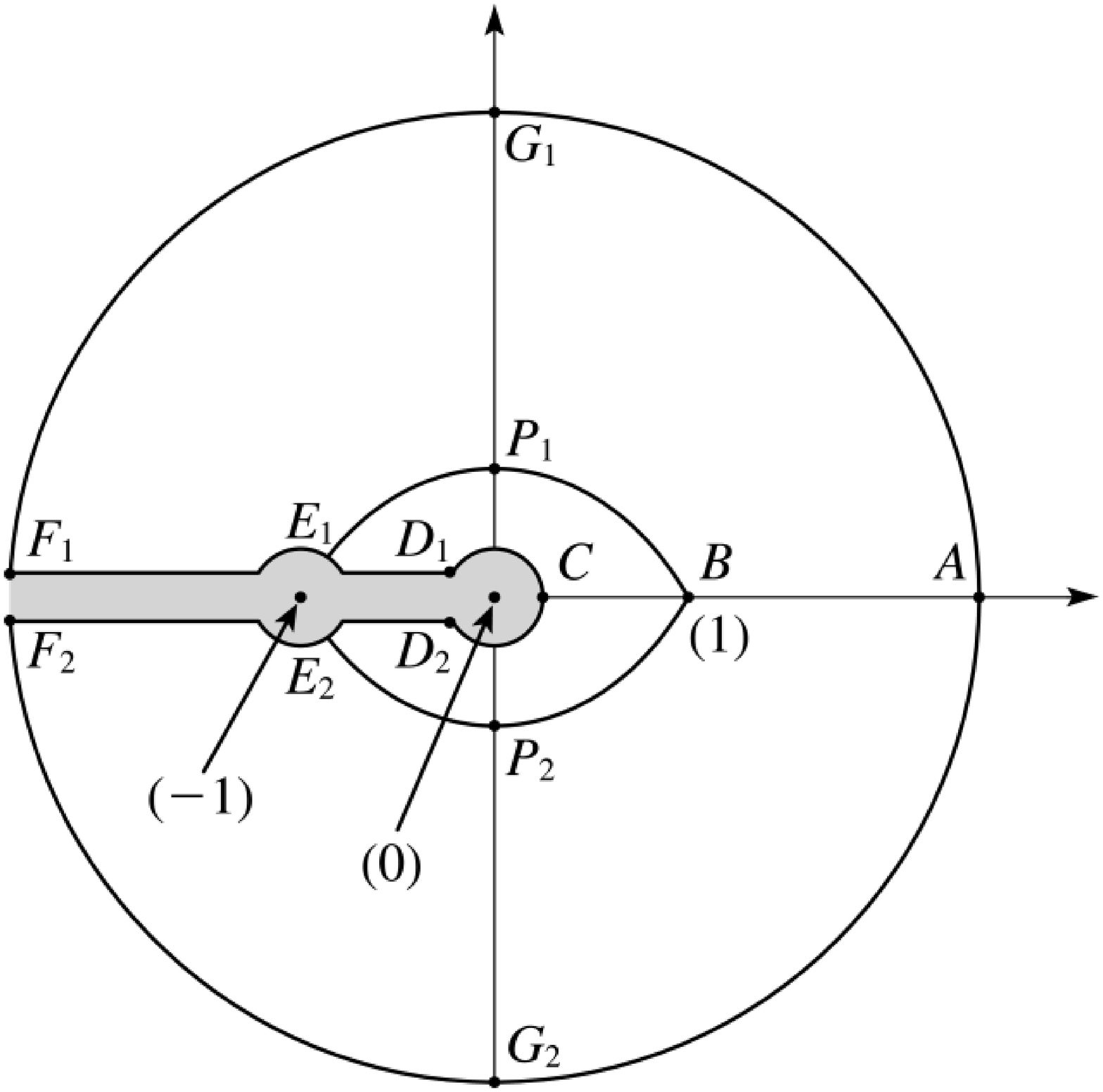,angle=0,height=5.5cm,width=5.5cm}}}
\end{minipage}
\hspace*{3cm} 
\begin{minipage}{3cm}
\centerline{\protect\hbox{\psfig{file=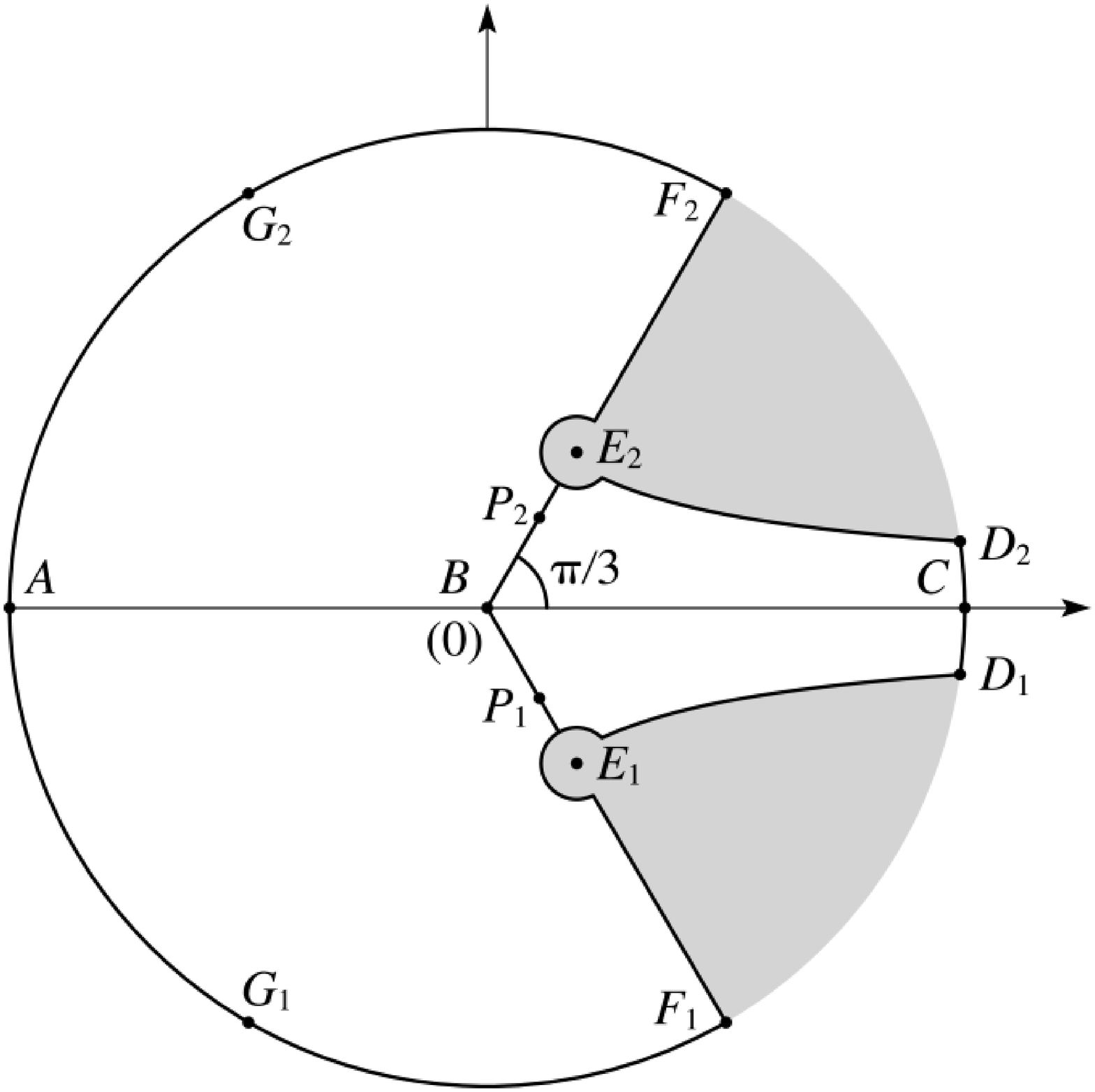,angle=0,height=5.5cm,width=5.5cm}}}
\end{minipage}
\end{center}
\caption{$z$-domain (left) and $\protect\zeta$-domain (right), with
corresponding points. These figures were taken from
http://dlmf.nist.gov/10.20; they are copyrighted by NIST and used with
permission.}
\label{fig1}
\end{figure}

From (\ref{eq2}), and taking the negative sign, let 
\begin{equation}
\frac{2}{3}\zeta ^{3/2}=\ln \left\{ \frac{1+\left( {1-z^{2}}\right) ^{1/2}}{z}\right\} 
-\left( {1-z^{2}}\right) ^{1/2},  \label{eq33}
\end{equation}
and 
$$
W=\zeta ^{-1/4}\left( \frac{1-z^{2}}{z^{2}}\right) ^{1/4}w.  
$$
The transformed variable $\zeta $ is real for real $z\in (0,1)$ ($\zeta \in
(0,+\infty )$), and $\zeta (z)$ can be defined by analytic continuation in
the whole complex plane cut along the negative real axis. This
transformation, and its correspondence in the $z$-plane, is depicted in
Fig. \ref{fig1}.

Observe that, if we are assuming that the principal values are taken, the
expression (\ref{eq33}) should not be used in all the complex plane and, for
instance, for real $z>1$ we would have complex values of $\zeta $, when $%
\zeta $ should be real for real $z>0$. This problem is solved by taking the
following: 
\begin{equation}
\frac{2}{3}(-\zeta )^{3/2}=\left(z^{2}-1\right)^{1/2} -\arctan\left\{\left(
z^{2}-1\right)^{1/2}\right\} ,  \label{eq34a}
\end{equation}
which in fact is an alternative when $\Re z>1$ (and in fact a larger
region). We can also write this in terms of logarithms as

\begin{equation}
\zeta (z)=\left\{ 
\begin{array}{l}
\left(\Frac{3}{2}\right)^{2/3} \left[ \log 
\Frac{1+\sqrt{1-z^{2}}}{z}-\sqrt{1-z^{2}}\right] ^{2/3},\,
\Re (z)\leq 1,z\notin {\mathbb{R}}^{-}, \\ 
-\left(\Frac{3}{2}\right)^{2/3} \left[ i\log 
\Frac{1+i\sqrt{z^{2}-1}}{z}+\sqrt{z^{2}-1}\right] ^{2/3},\,\Re (z)>1 .
\end{array}
\right.  \label{eq34b}
\end{equation}

We get (\ref{eq3}), with $u$ replaced by $\nu $, and 
$$
\psi(\zeta) =\Frac{5}{16\zeta ^{2}} +\Frac{\zeta z^{2}(z^{2}+4)}{4 (z^{2}-1)^3}.
$$

From \cite[Chap. 11, \S 10]{Olver:1997:ASF}) we identify the asymptotic
solutions by 
\begin{equation}
J_{\nu }(\nu z)=\frac{1}{\nu ^{1/3}}\left( \frac{4\zeta }{1-z^{2}}\right)
^{1/4}W_{0}(\nu ,\zeta ),  \label{eq36}
\end{equation}%
\begin{equation}
H_{\nu }^{(1)}(\nu z)=\frac{2e^{-\pi i/3}}{\nu ^{1/3}}\left( \frac{4\zeta }
{1-z^{2}}\right) ^{1/4}W_{-1}(\nu ,\zeta ),  \label{eq37}
\end{equation}
and 
\begin{equation}
H_{\nu }^{(2)}(\nu z)=\frac{2e^{\pi i/3}}{\nu ^{1/3}}\left( \frac{4\zeta }
{1-z^{2}}\right) ^{1/4}W_{1}(\nu ,\zeta ). 
\end{equation}

For the corresponding Liouville-Green expansions, set 
$$
\xi =\frac{2}{3}\zeta ^{3/2},  
$$
and define 
$$
V=\left( \frac{1-z^{2}}{z^{2}}\right) ^{1/4}w.  
$$

For our purposes we only need to consider $\left\vert\arg(z) \right\vert
\leq \frac{1}{2}\pi$, and so $\arg (\zeta) \leq 0$. Thus, taking $\arg (\xi)
=0$ when $\arg (\zeta)=0$, and letting $\xi$ depend continuously on $\zeta$
in the region corresponding to 
$\left\vert\arg(z)\right\vert\leq \frac{1}{2}\pi $, we find that $\arg (\xi) \leq 0$.

We observe then that $\xi \rightarrow +\infty $ as $z\rightarrow 0^{+}$ 
($\zeta \rightarrow +\infty )$. Also, $\left\vert \xi \right\vert \rightarrow
\infty $ as $z\rightarrow i\infty $, such that 
$\xi =iz- \frac{1}{2} i\pi +\mathcal{O}\left(z^{-1}\right)$ 
($\arg (\xi) \rightarrow -\pi )$. In
addition, $\left\vert \xi \right\vert \rightarrow \infty $ as $z\rightarrow
-i\infty $, such that $\xi =iz+ \frac{1}{2} i\pi +\mathcal{O}\left(
z^{-1}\right)$ ($\arg ( \xi ) \rightarrow -2\pi )$.

Then we have the transformed equation (\ref{eq3}), again with $u$ replaced
by $\nu $, and 
$$
\phi (\xi )=\frac{z^{2}(z^{2}+4)}{4(z^{2}-1)^{3}},  
$$

Now consider the coefficients given by (\ref{eq12}) - (\ref{eq14}). In
general, we prefer to numerically evaluate with respect to $z$, since $\xi $
is given explicitly in terms of $z$, but not vice versa. Thus, on account of 
$$
\frac{dz}{d\xi }=-\frac{z}{(1-z^{2})^{1/2}},  
$$
we have $E_{s}(\xi )=\hat{E}_{s}(z)$ (say), where 
\begin{equation}
\hat{E}_{s}(z)=\int_{z}^{\infty }t^{-1}(1-t^{2})^{1/2}\hat{F}_{s}(t)dt\quad
(s=1,2,3,\cdots ),  \label{eq43}
\end{equation}
in which 
$$
\hat{F}_{s}(z)=F_{s}\left( \xi (z)\right) \quad (s=1,2,3,\cdots ).
$$

Thus 
\begin{equation}
\hat{F}_{1}(z)=\frac{z^{2}(z^{2}+4)}{8(z^{2}-1)^{3}},\, 
\hat{F}_{2}(z)=\frac{z}{2( 1-z^{2})^{1/2}}\hat{F}_{1}^{\prime }(z),
\end{equation}
and 
\begin{equation}  \label{recuFs}
\hat{F}_{s+1}(z)=\frac{z}{2( 1-z^{2}) ^{1/2}} 
\hat{F}_{s}^{\prime }(z)-\frac{1}{2}\sum_{j=1}^{s-1} 
\hat{F}_{j}(z)\hat{F}_{s-j}(z)\quad (s=2,3,\cdots) .
\end{equation}
For convenience we have chosen the integration constants so that all
coefficients vanish at $z=\infty $.

In the general case, the coefficients $\hat{E}_{s}( z)$ can be computed
numerically by quadrature via (\ref{eq43}). But for Bessel's equation it
turns out that the coefficients can be explicitly computed, and in
particular they have the expression 
\begin{equation}
\hat{E}_{s}(z)=\frac{P_{s}(z^{2})}{(1-z^{2})^{3s/2}},  \label{coebe}
\end{equation}
where $P_{s}(x)$ are polynomials of degree $s$ in $x$.

Before establishing this, we note for the odd terms that
$$
\hat{E}_{2j+1}(z)=\frac{1}{(1-z)^{1/2}} \left[\frac{P_{2j+1}(z^{2})}
{(1-z^{2})^{3j+1}(1+z)^{1/2}}\right] \ (j=0,1,2\cdots) ,
$$
where the term in the square brackets is meromorphic at $z=1$. Hence, on
expanding around the turning point $\zeta =0$ ($z=1$), we deduce that 
(\ref{eq30a}) holds with $\alpha_{2j+1}=0$, as desired.

In order to establish (\ref{coebe}), we first show that the coefficients 
$\hat{F}_{s}(z)$ have the form 
\begin{equation}
\hat{F}_{s}(z)=\frac{z^{2}Q_{s}(z^{2})}{(1-z^{2})^{\frac{3}{2}(s+1)}},
\label{fq}
\end{equation}
with $Q_{s}(t)$ a polynomial of degree not larger than $s$ in $t$.

We firstly observe that $\hat{F}_{1}(z)$ has this form. Next, we substitute 
(\ref{fq}) in (\ref{recuFs}) and then, for (\ref{fq}) and (\ref{recuFs}) to
hold, we have that the polynomials must satisfy 
\begin{equation}
Q_{s+1}(t)=\left[ 1+\frac{1}{2}(3s+1)t\right] Q_{s}(t)+t(1-t)Q_{s}^{\prime
}(t)-\frac{t}{2}\sum_{j=1}^{s-1}Q_{j}(t)Q_{s-j}(t),  \label{Qs}
\end{equation}
which, starting from $Q_{1}(t)=(1/2+t/8)$, shows that $Q_{n}$ is of degree
not larger than $n$ and that (\ref{fq}) holds with polynomials given by (\ref{Qs}).

Now, the coefficients $\hat{E}_{s}$ can be computed from 
$$
\hat{E}_{s}(z)=-\int_{a}^{z}\hat{F}_{s}(r)\frac{\left(1-r^{2}\right)^{1/2}}
{r}dr= -\int_{a}^{z}\frac{rQ_{s}(r^{2})}{(1-r^{2})^{1+3s/2}}dr,
$$
where the starting point $a$ is chosen depending on the solution to be
matched. These integrals can be explicitly computed. We have 
$$
\int \frac{rQ_{s}(r^{2})}{(1-r^{2})^{1+3s/2}}dr=
\frac{1}{2}\int \frac{{Q_{s}(t)}}{{(1-t)^{1+3s/2}}}dt,
$$
and because $Q_{s}$ is a polynomial of degree not larger than $s$ we are
left with integrals of the form 
\begin{equation}
\int \frac{t^{k}}{(1-t)^{a}}dt=H(t)+C,\,H(t)=\frac{R_{k}(t)}{(1-t)^{a-1}},
\label{primi}
\end{equation}
with $k<a-1$ and $R_{k}$ a polynomial of degree not larger than $k$ which
can be computed by differentiating (\ref{primi}). Because $k<a-1$, we have 
$H(\infty )=0$.

The polynomials $P_s$ in (\ref{coebe}) have the properties: 
$$
P_{2s}(0)=0,\,P_{2s+1}(0)=C_{2s+1},
$$
where $C_{2s+1}$ are the coefficients in the Stirling asymptotic series 
$$
\log {\rm \Gamma} (u)=\left( u-\frac{1}{2}\right) \log (u)-u+\frac{1}{2}\log
(2\pi )+\sum_{j=0}^{\infty }\frac{C_{2j+1}}{u^{2j+1}},\,u\rightarrow
\infty .
$$
We found that the property $P_{2s+1}(0)=C_{2s+1}$ holds by comparing 
(\ref{eq36}) with a similar expansion for $J_{\nu}(\nu z)$ but matching
the solution at $z=0$, using that $J_{\nu}(\nu z)\sim (\nu z /2)^{\nu}/{\rm \Gamma} (\nu+1)$
as $z\rightarrow 0$.

The first few polynomials are 
$$
\begin{array}{l}
P_{1}(x)=\frac{1}{24}(2+3x), \\ 
\\ 
P_{2}(x)=-\frac{1}{16}x(x+4), \\ 
\\ 
P_{3}(x)=\frac{1}{5760}(-16+1512x+3654x^{2}+375x^{3}), \\ 
\\ 
P_{4}(x)=\frac{1}{128}x(32+288x+232x^{2}+13x^{3}).
\end{array}
$$

Next, using 
$$
H_{\nu }^{(1)}(\nu z)\sim \left( \frac{2}{\pi \nu z}\right) ^{1/2}\exp
\left\{ i\nu z-\frac{1}{2}\nu \pi i-\frac{1}{4}\pi i\right\} ,  
$$
as $z\rightarrow \infty $ we match the solutions that are recessive in the
upper half $z$-plane, yielding 
\begin{equation}
H_{\nu }^{(1)}(\nu z)=-i\left( \frac{2}{\pi \nu }\right) ^{1/2}
\frac{1}{(1-z^{2})^{1/4}}V_{-1}(\nu ,\xi ),  \label{eq46}
\end{equation}
as $\nu \rightarrow \infty $, uniformly for $z$ lying in a domain which
contains the first quadrant (excluding a neighborhood of the turning point 
$z=1$).

Furthermore, we know 
$$
J_{\nu }(\nu z)=\frac{c_{0}(\nu )}{( 1-z^{2}) ^{1/4}}V_{0}(\nu ,\xi),  
$$
as $\nu \rightarrow \infty $ in a domain which contains 
$\mathrm{Re}(z)\geq 0 $, excluding a neighborhood of the interval $[1,\infty )$. Using 
$\xi =iz-\frac{1}{2}i\pi +\mathcal{O}\left( z^{-1}\right)$, (\ref{eq16}) and 
$$
J_{\nu }(\nu z)\sim -\left( \frac{1}{2\pi \nu z}\right)^{1/2}
\exp\left\{-i\nu z+\frac{1}{2}\nu \pi i+\frac{1}{4}\pi i\right\} ,  
$$
as $z\rightarrow i\infty $, we find that $c_{0}(\nu )=(2\pi \nu )^{-1/2}$,
and hence 
\begin{equation}
J_{\nu }(\nu z)=\left( \frac{1}{2\pi \nu }\right) ^{1/2}
\frac{1}{\left(1-z^{2}\right) ^{1/4}}V_{0}(\nu ,\xi ).  \label{eq49}
\end{equation}

We remark that the expansions (\ref{eq46}) and (\ref{eq49}) are a
reformulation of the Debye expansions \cite[\S 10.19(ii)]{Olver:2010:BF} for 
$z\in \left( 0,1\right) $, which also holds for certain complex values.
Debye expansions valid for $1<z<\infty $ are also given in this reference.

\section{Bessel's equation: turning point coefficient functions}

We now define the (scaled) turning point coefficient functions for Bessel's
equation. Taking $w_{-1}(u,z) =H_{\nu }^{(1)}(z)$ and 
$w_{1}(u,z) = H_{\nu}^{(2)}(z)$ in (\ref{31d}) these are defined by 
\begin{equation}
e^{\pi i/3}H_{\nu }^{(1)}(\nu z) = \mathrm{Ai}_{-1}\left(\nu
^{2/3}\zeta\right) \mathcal{A}(\nu ,z) +
\mathrm{Ai}_{-1}^{\prime }\left(\nu^{2/3}\zeta\right)\mathcal{B}(\nu ,z),  \label{eq50}
\end{equation}
and 
\begin{equation}
e^{-\pi i/3}H_{\nu }^{(2)}(\nu z) = \mathrm{Ai}_{1}\left(\nu^{2/3}\zeta\right) 
\mathcal{A}(\nu ,z)+ \mathrm{Ai}_{1}^{\prime }\left(\nu^{2/3}\zeta\right)\mathcal{B}(\nu ,z).  \label{eq50a}
\end{equation}
Comparing these two representations with (\ref{eq36}) and (\ref{eq37}) we
perceive that in (\ref{scaled}) the scaling function here is given by
\begin{equation}
\phi(\nu ,z) =\frac{2^{3/2}}{\nu ^{1/3}}\left(\frac{\zeta }{1-z^{2}}
\right)^{1/4} .  \label{phi}
\end{equation}

From (\ref{eq17}) and (\ref{eq49}) we also have 
\begin{equation}
2J_{\nu }(\nu z) =\mathrm{Ai}\left(\nu ^{2/3}\zeta\right)\mathcal{A}(\nu ,z)
+\mathrm{Ai}^{\prime }\left(\nu ^{2/3}\zeta\right)\mathcal{B}(\nu ,z).
\label{eq52}
\end{equation}
A similar relation for $-2Y_{\nu }(\nu z)$ is obtained by replacing 
$\mathrm{Ai}$ and its derivative by $\mathrm{Bi}$ and its derivative in (\ref{eq52}).

As described earlier, the idea for computing the slowly varying coefficients 
$\mathcal{A}(\nu ,z)$ and $\mathcal{B}{(\nu ,z)}$ in a region containing
the turning point $z=1$ 
\footnote{Around the turning point $z=-1$, and in general for $\Re (z)<0$ we can use
the continuation formulas in Sect. 10.11 of \cite{Olver:2010:BF}} is to
invoke Cauchy's integral formula (\ref{eq31}). In this, the integration is
taken along a positively oriented closed loop containing $z$ in its
interior, but not the origin; additionally, it should not cross the negative
real axis in the $z$-plane; in this way, we can guarantee that in the $\zeta$
plane we stay away from the shaded region in Fig. \ref{fig1} (right) and
therefore that all the functions appearing in the integration are analytic
in a domain containing the path of integration and its interior (and
therefore Cauchy's integral formula holds for this integration path). For
our purposes it will be sufficient to consider a circular path of
integration in (\ref{31a}), which can be centered or not at $z=1$, but which
in any case must contain $z=1$, but not $z=0$.

\label{LGC}

From\ (\ref{eq30}), (\ref{scaled}) and (\ref{phi}), we have the following
expansion which we use to compute $\mathcal{B}(\nu ,z)$ on the path of
integration

$$
\mathcal{B}(\nu ,z)\sim \Frac{2\sqrt{2}\exp\left\{\alpha(\nu ,z)\right\}}
{\nu ^{2/3}\zeta ^{1/4}(1-z^{2})^{1/4}}
\sinh\left\{\Frac{1}{\nu}\beta(\nu,z)\right\} ,  
$$
where, formally,

$$
\alpha(\nu ,z)=\sum_{j=1}^{\infty }\frac{\hat{E}_{2j}(z)+d_{2j}(\xi )}{\nu^{2j}},  
$$
and 
$$
\beta (\nu ,z)=\sum_{j=0}^{\infty }\frac{\hat{E}_{2j+1}(z)-d_{2j+1}(\xi )}{\nu ^{2j}},  
$$
in which 
$$
d_{s}(\xi )=a_{s}/(s\xi ^{s}).  
$$
where $a_{1}=a_{2}=\frac{5}{72},$ and subsequent terms given by (\ref{A7}).

We observe that the cancellation of the coefficient as $\nu \rightarrow
\infty $ is located in the $\sinh $ term. In order to avoid loss of
accuracy, it is better to write 
\begin{equation}
\mathcal{B}(\nu ,z)\sim \Frac{2\sqrt{2}\beta (\nu ,z)\exp\left\{ \alpha (\nu ,z)\right\}}
{\nu^{5/3}\zeta^{1/4}(1-z^{2})^{1/4}}\mathrm{sinhc}\left\{\Frac{\beta(\nu,z)}{\nu}\right\} ,  
\label{Bsinhc}
\end{equation}
where the function $\mathrm{sinhc}(x)=\mathrm{sinh}(x)/x$ can be computed
for small $x$ with the Maclaurin series 
$$
\mathrm{sinhc}(x)=\sum_{k=0}^{\infty }\frac{x^{2k}}{(2k+1)!}.  
$$

Several properties that could be expected for the coefficient 
$\mathcal{B}(\nu ,z)$ can be seen to hold explicitly from (\ref{Bsinhc}):

\begin{enumerate}
\item {It is real for real values of $z$.}

\item {}$\mathcal{B}(\nu ,z)=\mathcal{O}(\nu ^{-5/3})$ as $\nu \rightarrow\infty $, 
more specifically

$$
\begin{array}{ll}
\mathcal{B}(\nu ,z) & =\Frac{2\sqrt{2}}{\nu^{5/3}\zeta ^{1/4}(1-z^2)^{1/4}}
(\hat{E}_{1}(z)-d_{1}(\xi ))F(\nu ,z) \\ 
& =\Frac{2\sqrt{2}}{\nu^{5/3}\zeta^{1/4}(1-z^{2})^{1/4}}
\left(\Frac{2+3z^{2}}{24(1-z^{2})^{3/2}}-\Frac{5}{72\xi }\right)F(\nu ,z),
\end{array}
$$
where $F(\nu ,z)=1+\mathcal{O}(\nu ^{-2})$.

\item {Only even powers of $\nu ^{-1}$ appear} if $F(\nu ,z)$ is expanded in
powers of $\nu ^{-1}$.
\end{enumerate}

Proceeding similarly with the $\mathcal{A}(\nu,z)$ coefficient, we get

\begin{equation}
\mathcal{A}(\nu ,z)\sim \Frac{2\sqrt{2}\zeta ^{1/4}
\exp\left\{ \tilde{\alpha}(\nu ,z)\right\}}
{\nu^{1/3}(1-z^{2})^{1/4}}\cosh\left\{\Frac{1}{\nu}\tilde{\beta}(\nu ,z)\right\} ,  \label{Acosh}
\end{equation}
where 
$$
\tilde{\alpha}(\nu ,z)=\sum_{j=1}^{\infty }\frac{\hat{E}_{2j}(z)
+\tilde{d}_{2j}(\xi )}{\nu ^{2j}},  
$$

$$
\tilde{\beta}(\nu ,z)=\sum_{j=0}^{\infty }\frac{\hat{E}_{2j+1}(z)-
\tilde{d}_{2j+1}(\xi )}{\nu ^{2j}},  
$$
and 
$$
\tilde{d}_{s}(\xi )=\tilde{a}_{s}/(s\xi ^{s}),  
$$
in which $\tilde{a}_{1}=\tilde{a}_{2}=-\frac{7}{72}$, and subsequent terms
satisfying the same recursion formula (\ref{A7}) as for the $a_{s}$
coefficients.

For the function $\mathcal{A}(\nu ,z)$, we have the expected properties from
(\ref{Acosh}):

\begin{enumerate}
\item {It is real for real values of $z$.}

\item {}$\mathcal{A}(\nu ,z)=\mathcal{O}(\nu ^{-1/3})$ as $\nu \rightarrow
\infty $, more specifically 
\begin{equation}
\mathcal{A}(\nu ,z)=\frac{2\sqrt{2}\zeta ^{1/4}G(\nu ,z)}{\nu^{1/3}(1-z^{2})^{1/4}},
\end{equation}
where $G(\nu ,z)=1+\mathcal{O}(\nu ^{-2})$.

\item {Only even powers of $\nu ^{-1}$ appear} if $G(\nu,z)$ is expanded in
powers of $\nu ^{-1}$.
\end{enumerate}

For the derivatives ${H_{\nu }^{(1)}}^{\prime }(\nu z)$ and 
${H_{\nu }^{(2)}}^{\prime }(\nu z)$ we proceed as before and we express them in the form 
\begin{equation}
e^{i\pi /3}\nu {H_{\nu }^{(1)}}^{\prime }(\nu z)=\mathrm{Ai}_{-1}(\nu
^{2/3}\zeta )\mathcal{C}(\nu ,z)+\mathrm{Ai}_{-1}^{\prime }(\nu ^{2/3}\zeta )
\mathcal{D}(\nu ,z),  \label{definco3}
\end{equation}
and 
\begin{equation}
e^{-i\pi /3}\nu {H_{\nu }^{(2)}}^{\prime }(\nu z)=\mathrm{Ai}_{1}(\nu
^{2/3}\zeta )\mathcal{C}(\nu ,z)+\mathrm{Ai}_{1}^{\prime }(\nu ^{2/3}\zeta )
\mathcal{D}(\nu ,z).  \label{definco4}
\end{equation}

Solving (\ref{definco3}) and (\ref{definco4}) for $\mathcal{C}(\nu ,z)$ and 
$\mathcal{D}(\nu ,z)$, and considering the Airy-type expansions for the
derivatives \cite[10.20.9]{Olver:2010:BF}, we find that 
$$
\mathcal{C}(\nu ,z)=\mathcal{O}(\nu ^{-1/3}),\,\mathcal{D}(\nu ,z)=\mathcal{O}(\nu ^{1/3}),
$$
as $\nu \rightarrow \infty $. From the behavior noted above of 
$\mathcal{A}(\nu ,z)$ and $\mathcal{B}(\nu ,z)$ for large $\nu $ (which also hold for
their derivatives), we see that in the dominant coefficient $\mathcal{D}(\nu,z)$ 
the first term in (\ref{coeffd}) is the largest, which is indeed 
$\mathcal{O}(\nu ^{1/3})$, while for the coefficient $\mathcal{C}(\nu ,z)$
both terms in (\ref{coeffc}) are $\mathcal{O}(\nu ^{-1/3})$, which is the
correct order of this coefficient. Therefore, once the cancellation for the
coefficient $\mathcal{B}(\nu ,z)$ is avoided, no cancellations occur as $\nu
\rightarrow \infty $, and we expect both (\ref{coeffc}) and (\ref{coeffd})
to be numerically stable.

There are two possible approaches in numerically evaluating the coefficients 
$\mathcal{C}(\nu ,z)$ and $\mathcal{D}(\nu ,z)$. With $\mathcal{A}(\nu ,z)$
and $\mathcal{B}(\nu ,z)$ computed on the path of integration as described
above, the first method is to numerically evaluate the integrals (\ref{eq57}) 
and (\ref{eq58}).

Alternatively, one can use (\ref{coeffc}) and (\ref{coeffd}), with 
${\mathcal{A}}(\nu ,z)$ and $\mathcal{B}(\nu ,z)$, and their derivatives,
computed via Cauchy integral formulas. Then, $\zeta $ and $\zeta ^{\prime }$
can be evaluated directly from (\ref{eq33}), (\ref{eq34a}) and (\ref{eq34b}). 
There is, however, a possible loss of accuracy as $z\rightarrow 1$ in the
computation of $\zeta ^{\prime }(z)$ but it can be easily eliminated. We
have 
$$
\zeta ^{\prime }(z)=-\frac{\left( 1-z^{2}\right) ^{1/2}}{z\zeta ^{1/2}},
$$
and both the numerator and the denominator tend to $0$ as $z\rightarrow 1$
which implies loss of accuracy. In order to avoid this we put 
$\delta =\sqrt{1-z^{2}}$ and when $|\delta |$ is small we consider the 
Maclaurin series for 
$$
f(\delta )=\frac{3^{1/3}}{\delta }\left( \log \left( \frac{1+\delta }
{\sqrt{1-\delta ^{2}}}\right) -\delta \right)^{1/3}= 1+\frac{1}{5}d^{2}+
\frac{18}{175}d^{4}+\ldots  
$$
and compute 
\begin{equation}
\zeta ^{\prime }(z)=-\frac{2^{1/3}}{zf(\delta )}.  \label{eq61}
\end{equation}

Returning to $\mathcal{A}(\nu ,z)$ and $\mathcal{B}(\nu ,z)$, we end this
section describing a more direct, but less stable, method for their
computation (again, on the path of integration of (\ref{eq31})). 
From (\ref{eq50}) and (\ref{eq50a}) we have the exact representations

\begin{equation}
\begin{array}{ll}
\mathcal{A}(\nu ,z)=-2\pi i & \left\{ e^{i\pi /3}H_{\nu }^{(1)}(\nu z)
\mathrm{Ai}_{1}^{\prime }(\nu ^{2/3}\zeta )\right. \\ 
& \left. -e^{-i\pi /3}H_{\nu }^{(2)}(\nu z)
\mathrm{Ai}_{-1}^{\prime }(\nu^{2/3}\zeta )\right\} ,
\end{array}
\label{AHankel}
\end{equation}
and 
\begin{equation}
\begin{array}{ll}
\mathcal{B}(\nu ,z)=2\pi i & \left\{ e^{i\pi /3}H_{\nu }^{(1)}(\nu z)
\mathrm{Ai}_{1}(\nu ^{2/3}\zeta )\right. \\ 
& \left. -e^{-i\pi /3}H_{\nu }^{(2)}(\nu z)\mathrm{Ai}_{-1}(\nu ^{2/3}\zeta)\right\} .
\end{array}
\label{BHankel}
\end{equation}
Now, because we are assuming that a method to compute the Airy functions is
available (we are precisely considering expansions in terms of Airy
functions), we could compute these coefficients on the Cauchy contour
without the need to substitute the Airy functions by their Liouville-Green
expansions, as done before. For examples of Fortran implementations of
complex Airy functions see \cite{Amos:1986:A6A,Fabijonas:2004:A8A,Gil:2002:A8A}.

As before, the integration path is chosen in such a way that the functions
appearing in (\ref{AHankel}) and (\ref{BHankel}) can be computed without
recourse to the Airy-type expansion. With respect to the cylinder functions
involved in the computation of the coefficients, we use (\ref{eq46}) and 
(\ref{eq49}), along with $H_{\nu }^{(2)}(\nu z)=\overline{H_{\nu }^{(1)}(\nu 
\overline{z})}$.

In order to compute the coefficients in a numerically stable way along the
contour of integration, we need to verify that the two terms are not
suffering cancellations, which would happen if the two terms are
exponentially large and cancel each order. If both Hankel functions are
large, one should be replaced by the $J$ Bessel function. For example, we
see that inside the eye-shaped curve (Fig \ref{fig1}), the expression for 
$\mathcal{B}(\nu ,z)$ in (\ref{BHankel}) is unstable because two
exponentially large quantities are subtracting: all terms in this expression
are dominant inside the eye-shaped curve. Outside the eye, this expression
is in principle stable because in each term a dominant function multiples a
recessive function.

Inside the eye we can write a satisfactory expression using 
$H_{\nu}^{(2)}(\nu z)=2J_{\nu }(\nu z)-H_{\nu }^{(1)}(\nu z)$. We obtain 
\begin{equation}
\mathcal{B}(\nu ,z)=2\pi i\left\{ H_{\nu }^{(1)}(\nu z)
\mathrm{Ai}_{0}(\nu^{2/3}\zeta ) -2e^{-i\pi /3}J_{\nu }(\nu z) 
\mathrm{Ai} _{-1}(\nu^{2/3}\zeta)\right\} .  \label{secondB}
\end{equation}

This expression not only can be used inside the eye-shaped region, but also
for the rest of the half plane $\Im (z)\geq 0$; however, close to the real
axis when $z>1$ we get higher accuracy from Liouville-Green expansions using
(\ref{BHankel}). In our numerical algorithms we use (\ref{BHankel}) when 
$\Re \left( z\right) >1$ and (\ref{secondB}) in the rest of this half plane.
But, as described before, if we also expand the Airy functions instead of
computing them separately, we obtain asymptotic expansions for the
coefficient away from the turning point and switching from one expression to
the other for the coefficients is not needed in this case. See Sect. \ref{LGC}.

Similarly for $\mathcal{A}(\nu ,z)$ we have 
$$
\mathcal{A}(\nu ,z)=-2\pi i\left\{ H_{\nu }^{(1)}(\nu z)
\mathrm{Ai}_{0}^{\prime}(\nu ^{2/3}\zeta ) -2e^{-i\pi /3}J_{\nu }(\nu z) 
\mathrm{Ai}_{-1}^{\prime }(\nu ^{2/3}\zeta )\right\} .  
$$

We could do analogous substitutions to get a satisfactory formula when 
$\Im z\leq 0$. However this is not really necessary because the coefficients 
$\mathcal{A}(\nu ,z)$ and $\mathcal{B}(\nu ,z)$ are real on the real line
and, as commented before, analytic in a domain containing the integration
path. Therefore, by Schwarz reflection principle, in this domain 
$\mathcal{A}(\nu ,\bar{z})=\overline{\mathcal{A}(\nu ,z)}$ and 
$\mathcal{B}(\nu ,\bar{z})=\overline{\mathcal{B}(\nu ,z)}$. 
From a computational point of view, this
reduces by one half the complexity of evaluating the Cauchy integral,
provided we take a symmetric contour with respect to the real axis, as we
will do.

This more direct approach has some disadvantages with respect to the
expansions in terms of stability. First, we notice that, even when the
solutions are chosen adequately, there still remains some cancellation in
the $\mathcal{B}(\nu,z)$ coefficient for large orders. Additionally, there is also
some accuracy degradation in the computation of Airy functions for large
arguments due to unavoidable loss of accuracy in the computation of
exponentials of large argument, as we next describe.

\section{Numerical results}

We now give several numerical illustrations for the performance of both the
Liouville-Green expansions and the approximation around the turning point.
For this purpose, we have coded our algorithms in Fortran 90, and compared
with the values given by Amos' algorithm (which has typically 13-14 digits
accuracy).  Amos' algorithm uses a variety of methods for computing
Bessel functions, depending on the values of $\nu$ and $z$, and in particular
Airy-type expansions close to the turning point. 
We have also tested our methods using Maple$^{TM}$. 
The comparison with Amos' algorithm is
used as an exhaustive testbench of the numerical stability of our approximations
in fixed precision arithmetic; both Amos' program and our the implementation
are fast enough to provide many thousands of function values in just a second. The
tests in variable precision with Maple$^{TM}$ are much slower and much less
exhaustive, but they will allow us to explore higher accuracies.

\subsection{Testing of the new Liouville-Green approximations}

Here we numerical evaluate the new expansions (\ref{eq46}) and (\ref{eq49}).
We notice that the analytic continuation formulas of 
\cite[\S 10.11]{Olver:2010:BF} can be used to compute cylinder functions for $\Re z<0$ from
the values for $\Re z>0$. For instance, we have 
$$
J_{\nu }(ze^{\pm \pi i})=e^{\pm \nu \pi i}J_{\nu }(z).  
$$
Therefore, testing for $\Re z>0$ is enough. However, we are also considering
the case $\Re z<0$ for the case of $J_{\nu }(z)$ in order to show the full
validity region of the expansion (\ref{eq49}).

\begin{figure}[tbp]
\hspace*{-2cm} \epsfxsize=19cm \epsfbox{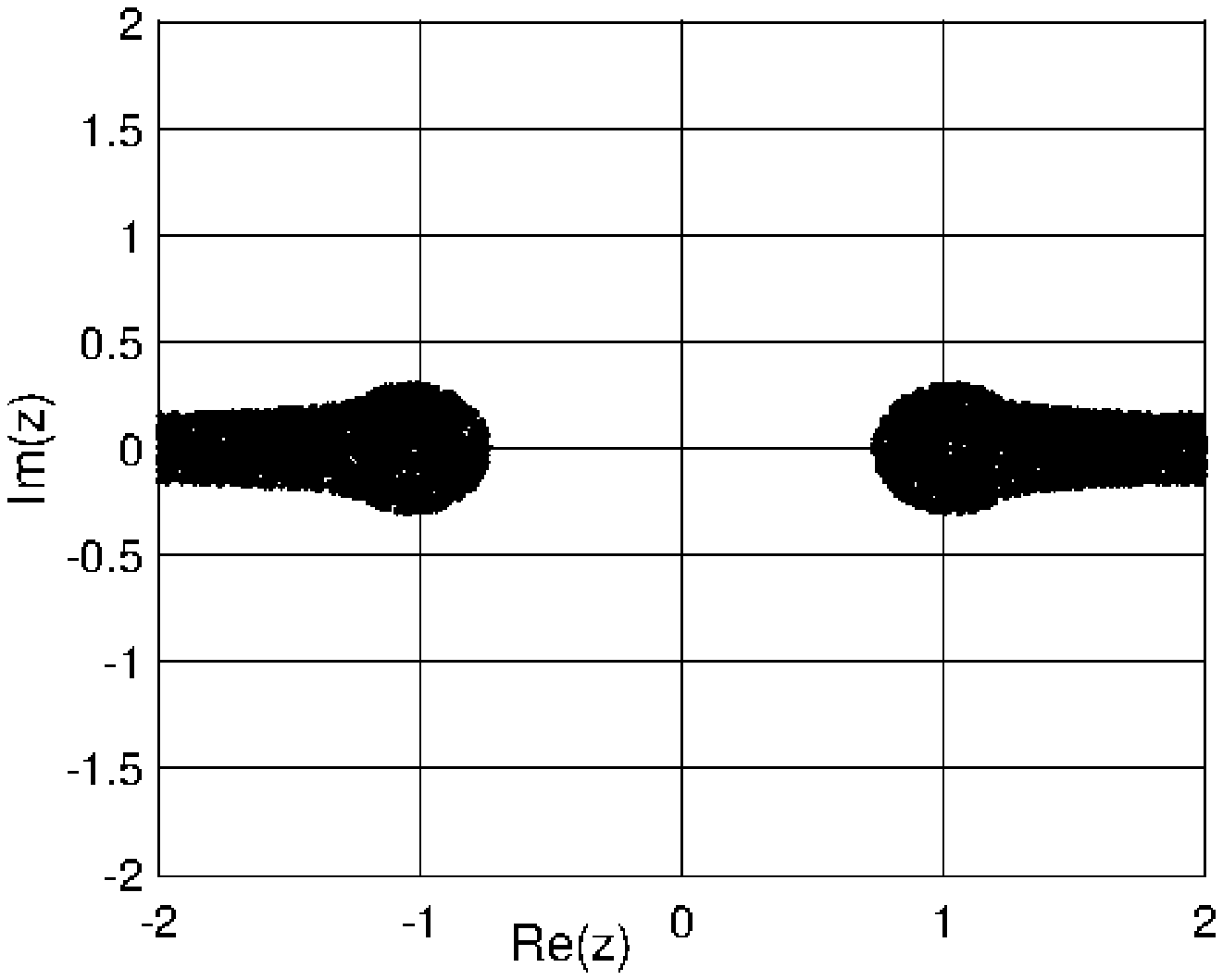}
\caption{ Comparison of the function values obtained with the expansion 
(\protect\ref{eq49}) against those obtained with Amos' algorithm 
\protect\cite{Amos:1986:A6A} for computing $J_{\protect\nu}(\protect\nu z)$ 
with $\protect\nu=100$ and $-2<\Re{z}<2,\,-2<\Im z<2$. The points where the 
relative error is greater than $10^{-12}$ are plotted.}
\label{fig:fig02}
\end{figure}

A test of the performance of the truncated expansion 
\begin{equation}
J_{\nu }(\nu z)\approx \left( \frac{1}{2\pi \nu }\right)^{1/2}
\frac{1}{\left(1-z^{2}\right)^{1/4}}
\exp\left\{ -\nu \xi +\displaystyle\sum_{s=1}^{n}(-1)^{s}
\frac{E_{s}(\xi )}{\nu ^{s}}\right\} ,  \label{Japprox}
\end{equation}
can be seen in Fig. \ref{fig:fig02}. The figure shows the comparison of
the function values obtained with $n=14$ in (\ref{Japprox}) against those
obtained with Amos' algorithm \cite{Amos:1986:A6A} for computing 
$J_{\nu}(\nu z)$ with $\nu =100$. Fig. \ref{fig:fig02} 
shows the comparison for
random values of the variable $z$ generated in the domain 
$-2<\Re {z}<2,\,-2<\Im z<2$, respectively. The points where the relative error is
greater than $10^{-12}$ are plotted in the figures. As the figure shows, and
as expected, the Liouville-Green expansion (\ref{eq49}) loses accuracy close to the turning
point and for real values of $z$ with $|z|>1$. It is worth noting than in 
the neighborhood of the turning point we can consider our expansions with 
coefficients computed via Cauchy integrals that we are discussing next, while
for $z>1$ we can compute $J_{\nu}(\nu z)$ using its relation with Hankel functions
and the Liouville-Green expansions for these functions or, alternatively, we can use
(\ref{eq52}) with coefficients computed from its asymptotic approximation. It appears 
then that for $\nu\ge 10$ it is possible to compute $J_{\nu}(\nu z)$ in the whole 
complex $z$-plane with around $15$ digits accuracy only by resorting to asymptotic 
approximations.

\begin{figure}[tbp]
\hspace*{-2cm} \epsfxsize=18cm \epsfbox{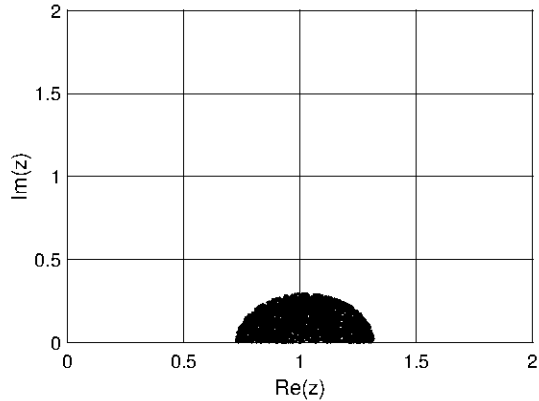}
\caption{ Comparison of the function values obtained with the expansion 
(\ref{eq46}) against those obtained with Amos' algorithm \cite{Amos:1986:A6A} 
for computing $H^{(1)}_{\nu}(\nu z)$ with $\nu=100$ and 
$0<\Re{z}<2,\,0<\Im z<2$. The points where the relative
error is greater than $10^{-12}$ are plotted.}
\label{fig:fig03}
\end{figure}

Finally, the performance of the Liouville-Green approximation 
\begin{equation}
H_{\nu }^{(1)}(\nu z) \approx -i\left(\frac{2}{\pi \nu}\right)^{1/2} 
\frac{1}{\left( {1-z^{2}}\right) ^{1/4}} \exp \left\{\nu \xi +
\displaystyle\sum_{s=1}^{n}\frac{E_{s}(\xi) }{\nu ^{s}}\right\} ,  \label{H1approx}
\end{equation}
for $\Im z>0$, with again $n=14$, is illustrated in Fig. \ref{fig:fig03}.
As expected, the expansion fails in the proximity of $z=1$.

\subsection{Airy-type expansions via Cauchy's integral formula}
\label{ATE}

We now test the accuracy in the computation of the Airy type expansion 
(\ref{eq50}) for $H_{\nu }^{(1)}(\nu z) $. We employ two different approaches in
the approximation of the coefficient functions $\mathcal{A}(\nu ,z)$ and 
$\mathcal{B}(\nu ,z)$, both of which use the Cauchy integral formulas 
(\ref{eq31}). In the first approach, on the Cauchy contour we use the
Liouville-Green expansions only for the cylinder functions, and we assume an
algorithm for computing complex Airy functions is available; we have used
both Amos' algorithm \cite{Amos:1986:A6A} and our own algorithm 
\cite{Gil:2002:A8A}, with similar results. Thus, in (\ref{eq50}) 
$\mathcal{A}(\nu,z)$ and $\mathcal{B}(\nu ,z)$ (expressed by the appropriate Airy and
cylinder functions) are approximated by using the Liouville-Green expansions
(\ref{Japprox}) and (\ref{H1approx}) in the integrands of (\ref{eq31}), but
with the Airy functions computed from the algorithms cited.

The second approach uses the Liouville-Green approximations for both the
Airy functions and the cylinder functions in approximating 
$\mathcal{A}(\nu,z)$ and $\mathcal{B}(\nu ,z)$ in (\ref{eq50}) 
(although again we do use Airy algorithms to compute 
$\mathrm{Ai}_{-1}\left(\nu ^{2/3}\zeta\right)$
and $\mathrm{Ai}_{-1}^{\prime }\left(\nu ^{2/3}\zeta\right)$ when computing 
$H_{\nu }^{(1)}(\nu z)$ with (\ref{eq63})). Hence in this case, from 
(\ref{eq50}), (\ref{Bsinhc}) and (\ref{Acosh}), 
we have for $z$ inside $\mathcal{L}$ 
\begin{equation}
e^{i\pi /3}H_{\nu}^{(1)}(\nu z) \approx \mathrm{Ai}_{-1}\left(\nu^{2/3}\zeta\right) 
\mathcal{A}_{m}(\nu ,z) +\mathrm{Ai}_{-1}^{\prime}\left(\nu ^{2/3}\zeta\right) 
\mathcal{B}_{m}(\nu ,z),  \label{eq63}
\end{equation}
where 
\begin{equation}
\mathcal{A}_{m}(\nu ,z)=\frac{1}{2\pi i} \oint\limits_{\mathcal{L}} 
\Frac{2\sqrt{2}\zeta ^{1/4}\exp \left\{ \tilde{\alpha}_{m}(\nu,t)\right\}}
{\nu ^{1/3}(1-t^{2})^{1/4}( t-z)} 
\cosh \left\{\Frac{1}{\nu}\tilde{\beta}_{m}(\nu,t)\right\}
dt,  \label{eq64}
\end{equation}
\begin{equation}
\mathcal{B}_{m}(\nu ,z)=\frac{1}{2\pi i} \oint\limits_{\mathcal{L}} 
\Frac{2\sqrt{2}\beta _{m}(\nu,t)\exp \left\{ \alpha_{m}(\nu,t)\right\}}
{\nu ^{5/3}\zeta(t)^{1/4}(1-t^{2})^{1/4}(t-z)} 
\mathrm{sinhc}\left\{\Frac{\beta _{m}(\nu,t)}{\nu}\right\} dt,  \label{eq65}
\end{equation}
in which 
$$
\tilde{\alpha}_{m}(\nu,t)=\sum_{j=1}^{m} \frac{\hat{E}_{2j}(t) + 
\tilde{d}_{2j}\left(\xi(t)\right)}{\nu ^{2j}},  
$$
$$
\tilde{\beta}_{m}(\nu,t)=\sum_{j=0}^{m-1} \frac{\hat{E}_{2j+1}(t) - 
\tilde{d}_{2j+1}\left(\xi(t)\right)}{\nu ^{2j}},  
$$
$$
\alpha _{m}(\nu,t)=\sum_{j=1}^{m} \frac{\hat{E}_{2j} (t) +
d_{2j}\left(\xi(t)\right)}{\nu ^{2j}},  
$$
and 
$$
\beta _{m}(\nu ,t)=\sum_{j=0}^{m-1} \frac{\hat{E}_{2j+1}(t) -
d_{2j+1}\left(\xi(t)\right)}{\nu ^{2j}}.  
$$

\begin{figure}[tbp]
\epsfxsize=13cm \epsfbox{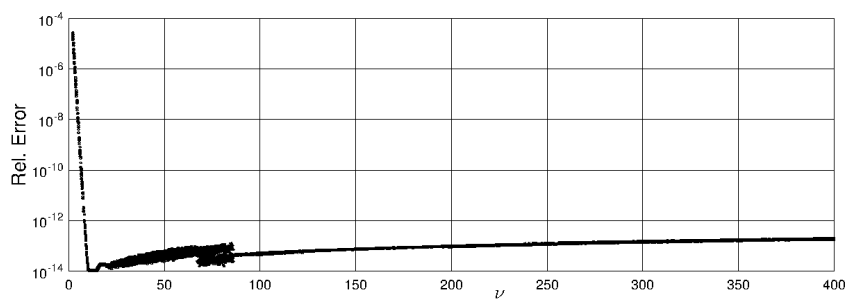}
\caption{ Relative errors in the comparison of the function values obtained
with the Airy-type expansion with $n$ terms for $H^{(1)}_{\nu}(\nu (1+0.1 i))$ against those obtained with Amos' algorithm \cite{Amos:1986:A6A}. In the evaluation of the coefficients the Airy
functions are also computed with \cite{Amos:1986:A6A}}
\label{fig:fig04}
\end{figure}

Two different tests of the accuracy of the Airy-type asymptotic expansion
for $H_{\nu }^{(1)}(\nu z)$ (with coefficients evaluated by Cauchy's
integral formula), for $z$ fixed or $\nu $ fixed, are considered. The
results are shown in Fig. \ref{fig:fig04}, for which the Airy function has
been evaluated both using \cite{Gil:2002:A8A} and \cite{Amos:1986:A6A}, with
similar results. In all cases we take $m=7$, which means we are considering
terms up to $\mathcal{O}(\nu^{-n})$, $n=2m=14$. Of course, similar tests can be
made for other solutions of the Bessel equation (for instance 
$J_{\nu}(\nu z)$ or $Y_{\nu}(\nu z)$), and the accuracy results are very similar,
as can be expected since all solutions are computed with the same coefficients.

In Fig. \ref{fig:fig04}, the relative error in the comparison of the
function values obtained with the Airy-type expansion for 
$H_{\nu}^{(1)}(\nu(1+0.1i)),\,2<\nu <400$ against those obtained with Amos' algorithm is
shown. Several different sources of error are apparent in the figure: the error
for $\nu $ small due to the Liouville-Green approximations for cylinder
functions used in the Airy-type expansions, the small increase in the
error due to rounding for $\nu $ large caused by the $\mathcal{B}(\nu ,z)$
coefficient, and the unavoidable loss of accuracy in the computation of
Airy functions for large arguments.

\begin{figure}[tbp]
\epsfxsize=13cm \epsfbox{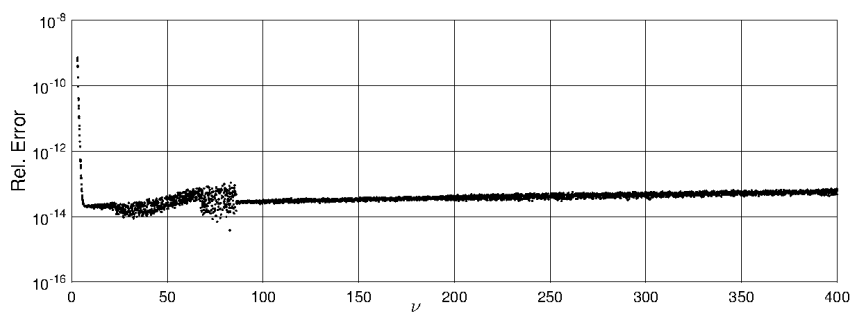}
\caption{ Relative errors in the comparison of the function values obtained
with the Airy-type expansion for $H^{(1)}_{\nu}(\nu (1+0.1 i))$ against those 
obtained with Amos' algorithm \cite{Amos:1986:A6A}. The asymptotic expansion 
of the coefficients is considered over the Cauchy contour. We take $n=14=2m$ 
in Eqs. (\ref{eq64}) and (\ref{eq65}).}
\label{fig:fig06}
\end{figure}

In Fig. \ref{fig:fig06} we show the same results as in Fig. \ref{fig:fig04},
but using the asymptotic expansions for the coefficients given in Sect. 
\ref{LGC}. We observe a clear improvement over Fig. \ref{fig:fig04}, particularly
for low $\nu$. As $\nu$ becomes larger, there a slow increase in accuracy
loss (smaller than in the previous case). This is due to the increasing
inaccuracies in the computation of Airy functions as the argument becomes
larger, which give rise to errors in the computation of the Hankel function
when using (\ref{eq50}). These inaccuracies are unavoidable in finite
precision arithmetic and, as described in \cite{Gil:2002:A8A}, can only be
removed by considering scaled functions (with the dominant exponential
factor exactly scaled out).

\begin{figure}[tbp]
\epsfxsize=13cm \epsfbox{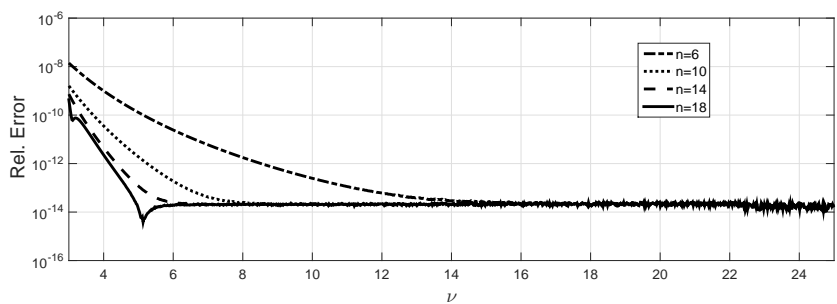}
\caption{ Same as Fig. \ref{fig:fig06} but for smaller $\nu$ and four selections of $n=2m$: $n=6,10,14,18$.}
\label{fig:fig07}
\end{figure}

\begin{figure}[tbp]
\epsfxsize=13cm \epsfbox{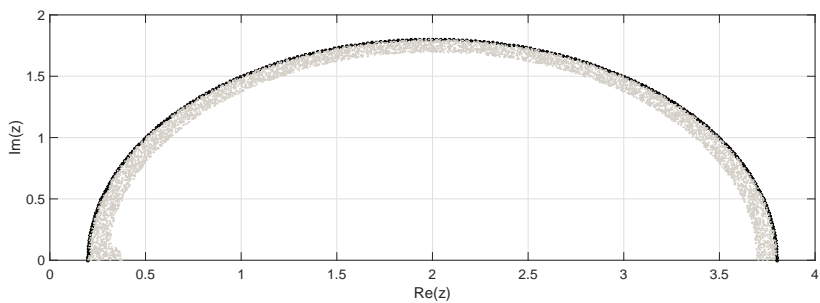}
\caption{ Comparison of the function values obtained with the Airy-type
expansion for $H^{(1)}_{10}(10z)$ against those obtained with Amos'
algorithm \cite{Amos:1986:A6A}. The points where the relative error
is greater than $10^{-13}$ are plotted. We take $n=2m=14$.}
\label{fig:fig05}
\end{figure}

Fig. \ref{fig:fig07} shows the same tests, but focusing in smaller values
of $\nu$ and for two selections in the number of terms for the expansion. We
observe two tendencies, particularly for $n=2m=18$: a decrease of the error
as $\nu$ increases due to the fact that the expansion becomes more accurate,
and a slow increase of the error due to the finite precision computation of
the Airy functions in the expression for $H^{(1)}_{\nu}(\nu z)$; this error increase
is not attributable to our approach. These tendencies give an optimal
accuracy for $\nu$ close to $5$. As commented, the accuracy for larger $\nu$
in finite precision arithmetic can be improved by considering scaled
functions.

In Fig. \ref{fig:fig05} we plot the points where the relative error is
greater than $10^{-13}$ in the comparison of the function values obtained
with the Airy-type expansion for $H^{(1)}_{10}(10 z)$ against those obtained
with Amos' algorithm. The random $z$ points in the test have been generated
inside the semi-circle limited by the integration contour used to compute
the coefficient functions $A(\nu,z)$ and $B(\nu,z)$ of the Airy-type
asymptotic expansion (a circle of center $z_c=2$ and radius $R=1.8$). As can be
seen in the figure, the accuracy obtained with the expansion is better than 
$10^{-13}$ in a large portion of the domain although, as expected, it worsens
when approaching the integration contour. This figure illustrates the
accuracy in the discretization of the Cauchy integral in a large region, but
with $z$ not too close to the contour of integration.

Preliminary tests show that the time spent for the computation of the Liouville-Green 
expansions in \S 3 and the Airy-type expansions in \S 4 is similar to the time for the implementations of Debye and Airy-type expansions in Amos' algorithm \cite{Amos:1986:A6A}. 
However, as expected, when the Cauchy's integral formula is used our approach becomes 
slower for computing a single function value, because we need to compute the coefficients 
a number of times over the contour.

The advantage of the Cauchy approach is that the most costly computation is the 
evaluation of the coefficients $E_{j}(\xi)$ and the numerators in the sums of 
Eqs. (\ref{eq29}) and (\ref{eq30}), but this computation is done
once and for all over the contour. 
Our approach permits the selection of a convenient Cauchy contour 
depending on the application. If many functions values are needed in 
a certain $z$-region, a good approach is to precompute the coefficients in a 
circuit containing this region (if it is possible and the shaded regions of Fig. \ref{fig1} 
can be avoided), and then applying the discretized version of Cauchy integral formula 
(Eq. (\ref{31c})) as many times as needed (and the recomputation
of the sums in Eqs. (\ref{eq29}) and (\ref{eq30}) if $\nu$ is also varying).

At this point, it is important to stress that the main interest of the
Cauchy technique lies in the computability more than in the efficiency (although it is efficient). 
It provides a new and direct method of computation of the coefficients
of Airy-type expansions with potential applications to more complicated cases.

\begin{figure}[tbp]
\epsfxsize=13cm \epsfbox{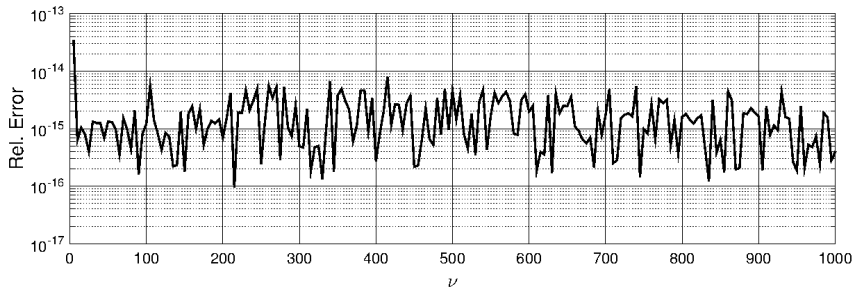}
\caption{ Error in the computation of $H_{\nu}^{(1)}(\nu z)$
for $z=1+0.1 i$ when the coefficients over the Cauchy contour are computed
with $16$ digits and using their LG asymptotic expansion. The rest of
computations are performed with $50$ digits using Maple. As before, $n=14$.}
\label{fig:fig08}
\end{figure}

In addition to the tests in fixed precision arithmetic, we have performed
additional test in variable precision using Maple$^{TM}$. Firstly, we have
checked that the coefficients in the Airy-type expansions are computed in a
numerically stable way with our scheme and that the remaining error
degradation in Fig. \ref{fig:fig06} is due to loss of accuracy in the
computation of the Airy functions when using (\ref{eq63}). For this purpose,
we have used Maple$^{TM}$ for computing the coefficients with a fixed number
of digits (we take 16 digits) and then computed the Hankel function using (%
\ref{eq63}) with a sufficiently high number of digits. When we test the
value of the Hankel function with, say $50$ digits, we observe that the
error is always close to $16$ digits accuracy, which shows that the
coefficients have been consistently computed with that accuracy and without
accuracy loss for high $\nu$. This is illustrated in Fig. \ref{fig:fig08}.

\begin{figure}[tbp]
\epsfxsize=13cm \epsfbox{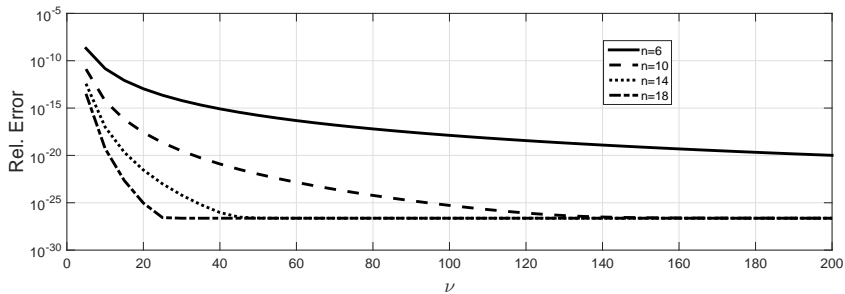}
\caption{ Error in the computation of $H_{\nu}^{(1)}(\nu z)$
for $z=1+0.1 i$ and various selections of $n=2m$ when the coefficients over
the Cauchy contour and the rest of computation are performed with $50$
digits using Maple. The limit in the minimal possible error (close to $%
10^{-26}$) is due to the discretization of the Cauchy integral with $N=500$
points}
\label{fig:fig09}
\end{figure}

In all the computations, we have used $500$ points in the contour of
integration. The error associated with the discretization of the Cauchy
integral is so small that it has no impact on the previous results. For
observing this error, a higher number of digits should be considered. Fig.
\ref{fig:fig09} shows the results when $500$ points over the Cauchy contour
are considered and the computations are done with $50$ digits.

As can be expected, because the integrand is periodic, the trapezoidal rule
has exponential convergence. Indeed, we have checked that smallest reachable
error roughly depends on the number of points over the Cauchy contour as 
$10^{-N/20}$.

We have checked numerically that when the effect of discretization of the
Cauchy integral is negligible, the relative error when $n$ terms in the
asymptotic expansion of the coefficients are considered, the relative error
in the computation of $H^{(1)}_{\nu} (\nu z)$ close to the turning point varies as 
$C_{n+1}/\nu^{n+1}$, where $C_{n+1}$ does not depend on $\nu$ and increases
with $n$. This is the expected behavior for an asymptotic expansion for
large $\nu$. Some estimations of these constants close to $z=1$ are shown in
table \ref{tabla1}

\begin{table}[tbp]
\[
\begin{array}{lllllllll}
\hline
n & 4 & 6 & 8 & 10 & 12 & 14 & 16 & 18 \\ \hline
C_{n+1} & 0.00015 & 0.00013 & 0.00021 & 0.00051 & 0.0018 & 0.0089 & 0.056 & 
0.46 \\ \hline
\end{array}
\]%
\caption{Computational error constants for the computation of $H^{(1)}_{\nu}(\nu z)$ 
close to $z=1$. The relative error is in good
approximation given by $C_{n+1}/\nu^{n+1}$}
\label{tabla1}
\end{table}

It would be important to be able to establish strict error bounds and to
compare them with these experimental errors. This should be achievable by
bounding the errors in the Liouville-Green expansions used for the
coefficients, similar to the error bounds given in \cite{Dunster:1998:AOT},
and again appealing to Cauchy's integral formula. This will be considered in
a subsequent paper.

\appendix
\section{Appendix: Exponential-form expansions for Airy functions}

{\normalsize
Firstly, the Airy functions $\mathrm{Ai}_{j}\left(u^{2/3}\zeta\right)$
satisfy
\begin{equation}
d^{2}w/d\zeta ^{2}=u^{2}\zeta w.  \label{A1}
\end{equation}

Letting $\xi =\int \zeta ^{1/2}d\zeta =\frac{2}{3}\zeta ^{3/2}$ 
(as in (\ref{eq7})), we then have that the functions 
$V=\zeta ^{1/4}\mathrm{Ai}_{j}\left( u^{2/3}\zeta \right) $ satisfy 
$$
\frac{d^{2}V}{d\xi ^{2}}=\left\{ u^{2}-\frac{5}{36\xi ^{2}}\right\} V.
$$

An asymptotic solution of the form (\ref{eq10}) now applies. In particular,
on identifying solutions of (\ref{A1}) that are recessive at $\zeta =+\infty$, 
we have that there exists a constant $c(u)$ such that 
$$
\mathrm{Ai}\left(u^{2/3}\zeta\right) \sim c(u) \zeta^{-1/4}
\exp \left\{-u\xi +\displaystyle\sum_{s=1}^{\infty } (-1)^{s}
\frac{e_{s}(\xi)}{u^{s}}\right\} ,  
$$
as $u^{2/3}\zeta \rightarrow \infty $ in the sector $\left\vert\arg
\left(\zeta \right)\right\vert \leq \pi -\delta $ ($\delta >0)$.

The coefficients in this expansion are given by (\ref{eq12}) - (\ref{eq14}),
with $E$ and $F$ replaced by $e$ and $f$, respectively, and 
$\phi (\xi )=-\frac{5}{{36}}\xi ^{-2}$. In particular, 
$e_{s}(\xi )=\int_{\infty }^{\xi}f_{s}(t)dt$, where 
$$
f_{1}(\xi )=-\frac{5}{72\xi ^{2}},\quad f_{2}(\xi )=-\frac{5}{72\xi ^{3}},
$$
and 
$$
f_{s+1}(\xi )=-\frac{1}{2}f_{s}^{\prime }(\xi )-\frac{1}{2}
\displaystyle\sum_{j=1}^{s-1}f_{j}(\xi )f_{s-j}(\xi )\quad (s\geq 2).  
$$

Thus 
$$
f_{s}(\xi) =-\frac{a_{s}}{\xi^{s+1}},\quad e_{s}(\xi) =\frac{a_{s}}{s\xi ^{s}}
\quad (s=1,2,3,\cdots) ,  
$$
where $a_{1}=a_{2}=\frac{5}{72}$, and 
\begin{equation}
a_{s+1}=\frac{1}{2}(s+1) a_{s} +\frac{1}{2}
\displaystyle\sum_{j=1}^{s-1}a_{j}a_{s-j}\quad (s\geq 2) .  \label{A7}
\end{equation}

From the well-known leading term
$$
\mathrm{Ai}\left(u^{2/3}\zeta\right) \sim \frac{e^{-u\xi }}
{2\pi^{1/2}u^{1/6}\zeta ^{1/4}} \quad \left(u^{2/3}\zeta \rightarrow\infty\right) ,  
$$
we obtain $c(u) =1/\left(2\pi ^{1/2}u^{1/6}\right)$, and hence we deduce
that 
\begin{equation}
\mathrm{Ai}\left(u^{2/3}\zeta\right) \sim \frac{1}{2\pi^{1/2}u^{1/6}\zeta^{1/4}}
\exp\left\{ -u\xi +\displaystyle\sum_{s=1}^{\infty}(-1)^s 
\frac{a_{s}}{su^{s}\xi ^{s}}\right\},  \label{A9}
\end{equation}
which is uniformly valid for $\left\vert\arg \left( \zeta \right)\right\vert
\leq \pi -\delta $ ($\delta >0$). Expansions for 
$\mathrm{Ai}_{\pm 1}\left(u^{2/3}\zeta\right)$ can be obtained directly from this.

Next, from differentiating (\ref{A1}) we find that 
$y=\zeta ^{-1/2}\mathrm{Ai}_{j}^{\prime } \left(u^{2/3}\zeta\right)$ satisfy 
$$
\frac{d^{2}y}{d\zeta ^{2}}=\left\{ u^{2}\zeta +\frac{3}{4\zeta ^{2}}\right\} y.  
$$
Thus, again with $\xi =\frac{2}{3}\zeta ^{3/2}$, we have that 
$\tilde{V}=\zeta ^{-1/4}\mathrm{Ai}_{j}^{\prime }\left( u^{2/3}\zeta \right) $ satisfy 
$$
\frac{d^{2}\tilde{V}}{d\xi ^{2}}=\left\{ u^{2}+\frac{7}{36\xi ^{2}}\right\} 
\tilde{V}.  
$$

Similarly to (\ref{A9}), on using \cite[(9.7.6)]{Olver:2010:AAR}, we deduce
that 
\begin{equation}
\mathrm{Ai}^{\prime }\left( u^{2/3}\zeta \right) \sim 
-\frac{u^{1/6}\zeta^{1/4}}{2\pi ^{1/2}}
\exp\left\{ -u\xi +\displaystyle\sum_{s=1}^{\infty}(-1)^{s}
\frac{\tilde{a}_{s}}{su^{s}\xi ^{s}}\right\} ,  \label{A12}
\end{equation}
for $\left\vert \arg (\zeta )\right\vert \leq \pi -\delta $,
where $\tilde{a}_{1}=\tilde{a}_{2}=-\frac{7}{72}$, and subsequent terms also
satisfy (\ref{A7}). Expansions for $\mathrm{Ai}_{\pm 1}^{\prime
}\left( u^{2/3}\zeta \right) $ can be obtained directly from this.
}

\begin{acknowledgements}
We thank the referees for a number of helpful suggestions.
\end{acknowledgements}

\bibliographystyle{spmpsci}
\bibliography{coebesRR}

\begin{thebibliography}{10}
\providecommand{\url}[1]{{#1}}
\providecommand{\urlprefix}{URL }
\expandafter\ifx\csname urlstyle\endcsname\relax
  \providecommand{\doi}[1]{DOI~\discretionary{}{}{}#1}\else
  \providecommand{\doi}{DOI~\discretionary{}{}{}\begingroup
  \urlstyle{rm}\Url}\fi

\bibitem{Amos:1986:A6A}
Amos, D.E.: Algorithm 644: a portable package for {B}essel functions of a
  complex argument and nonnegative order.
\newblock ACM Trans. Math. Software \textbf{12}(3), 265--273 (1986)

\bibitem{Amos:1977:C6S}
Amos, D.E., Daniel, S.L., Weston, M.K.: C{DC} 6600 subroutines {IBESS} and
  {JBESS} for {B}essel functions {$I_{\nu }(x)$} and {$J_{\nu }(x),$} {$x\geq
  0,$} {$\nu \geq 0$}.
\newblock ACM Trans. Math. Software \textbf{3}(1), 76--92 (1977)

\bibitem{Bornemann:2011:AAS}
Bornemann, F.: Accuracy and stability of computing high-order derivatives of
  analytic functions by {C}auchy integrals.
\newblock Found. Comput. Math. \textbf{11}(1), 1--63 (2011)

\bibitem{Boyd:1987:AEF}
Boyd, W.G.C.: Asymptotic expansions for the coefficient functions that arise in
  turning-point problems.
\newblock Proc. Roy. Soc. London Ser. A \textbf{410}(1838), 35--60 (1987)

\bibitem{Dunster:1998:AOT}
Dunster, T.M.: Asymptotics of the eigenvalues of the rotating harmonic
  oscillator.
\newblock J. Comput. Appl. Math. \textbf{93}(1), 45--73 (1998)

\bibitem{Fabijonas:2004:A8A}
Fabijonas, B.R.: Algorithm 838: {A}iry functions.
\newblock ACM Trans. Math. Software \textbf{30}(4), 491--501 (2004)

\bibitem{Gil:2002:A8A}
Gil, A., Segura, J., Temme, N.M.: Algorithm 819: {AIZ}, {BIZ}: two {F}ortran 77
  routines for the computation of complex {A}iry functions.
\newblock ACM Trans. Math. Software \textbf{28}(3), 325--336 (2002)

\bibitem{Gil:2007:NSF}
Gil, A., Segura, J., Temme, N.M.: Numerical methods for special functions.
\newblock SIAM, Philadelphia, PA (2007)

\bibitem{Olver:1997:ASF}
Olver, F.W.J.: Asymptotics and special functions.
\newblock AKP Classics. A K Peters Ltd., Wellesley, MA (1997).
\newblock Reprint of the 1974 original [Academic Press, New York]

\bibitem{Olver:2010:AAR}
Olver, F.W.J.: Airy and related functions.
\newblock In: N{IST} handbook of mathematical functions, pp. 193--213. U.S.
  Dept. Commerce, Washington, DC (2010)

\bibitem{Olver:2010:BF}
Olver, F.W.J., Maximon, L.C.: Bessel functions.
\newblock In: N{IST} handbook of mathematical functions, pp. 215--286. U.S.
  Dept. Commerce, Washington, DC (2010)

\bibitem{Temme:1997:NAF}
Temme, N.M.: Numerical algorithms for uniform {A}iry-type asymptotic
  expansions.
\newblock Numer. Algorithms \textbf{15}(2), 207--225 (1997)

\end{thebibliography}

\end{document}